\documentclass[11pt]{amsart}
\usepackage{amsfonts,amssymb}
\usepackage{amsmath}
\usepackage{amsthm,amscd,latexsym,mathrsfs,hyperref}

\newtheorem{theorem}{Theorem}[section]
\newtheorem{corollary}[theorem]{Corollary}
\newtheorem{lemma}[theorem]{Lemma}
\newtheorem{proposition}[theorem]{Proposition}

\newtheorem{assumption}{Assumption}

\newtheorem{theoremintro}{Theorem}

\theoremstyle{definition}
\newtheorem{definition}{Definition}[section]
\newtheorem{remark}{Remark}[section]

\DeclareMathOperator{\sat}{\mathrm{sat}}
\DeclareMathOperator{\hypo}{\mathrm{hypo}}

\DeclareMathOperator{\ran}{ran}
\DeclareMathOperator{\tr}{tr}
\DeclareMathOperator{\Lat}{Lat}

\begin{document}

\title[L\"owner's Theorem in several variables]{L\"owner's Theorem in several variables\thanks{This time for real.}}
\author[Mikl\'os P\'alfia]{Mikl\'os P\'alfia}
\address{Department of Mathematics, Sungkyunkwan University, Suwon 440-746, Korea.}
\address{Hungarian Academy of Sciences, 1051 Budapest, Sz\'echenyi Istv\'an sqr. 9, Hungary.}
\email{palfia.miklos@aut.bme.hu}



\date{\today}

\begin{abstract}
In this paper we establish a multivariable non-commutative generalization of L\"owner's classical theorem from 1934 characterizing operator monotone functions as real functions admitting analytic continuation mapping the upper complex half-plane into itself. The non-commutative several variable theorem proved here characterizes several variable operator monotone functions, not assumed to be free analytic or even continuous, as free functions that admit free analytic continuation mapping the upper operator poly-halfspace into the upper operator halfspace over an arbitrary Hilbert space. We establish a new abstract integral formula for them using non-commutative topology, matrix convexity and LMIs. The formula represents operator monotone and operator concave free functions as a conditional expectation of a Schur complement of a linear matrix pencil on a tensor product operator algebra. This formula is new even in the one variable case. The results can be applied to any of the various multivariable operator means that have been constructed in the last three decades or so, including the Karcher mean. Thus we obtain an explicit, closed formula for these operator means of several positive operators.
\end{abstract}

\maketitle

\section{Introduction}
In 1934 L\"owner proved his influential theorem on operator monotone functions of a real variable which states that all such functions are precisely given by holomorphic functions mapping the upper complex half-plane into itself \cite{lowner}. Many different proofs of this result appeared using different techniques \cite{bendatsherman,kraus,hansen3,hansenpedersen,sparr,wigner}. Some monographs comparing and giving details about the various techniques are \cite{bhatia,donoghue,schilling}. Most of the proofs are based on establishing an integral characterization formula for real functions that are operator monotone and then extending the domain of the function to the upper complex half plane by analytic continuation through the formula. These results were used in 1980 by Kubo and Ando \cite{kubo} to construct a two-variable theory of operator means, which in some sense provides a characterization for two-variable operator monotone functions with special properties. The article \cite{kubo} has proven to be an influential step towards developing a general theory of multivariable operator means, which are by definition, uniquely induced by a certain class of normalized multivariable operator monotone functions. In \cite{kubo} two-variable operator monotone functions with some additional properties are considered and characterized. In particular the two-variable geometric mean has received a great deal of attention, many authors were considering various different methods to extend the geometric mean to more than two non-commuting operator variables \cite{ando,AF,Ba08,bhatiaholbrook,BMP,lawsonlim0,lawsonlim,limpalfia,moakher,palfia,palfia2,palfia3}. The problem was also motivated by practical applications in medical imaging \cite{AF,moakher}, radar imaging technology \cite{Ba08} and many others. The various different techniques of extensions seemingly lead to different geometric means in more than two non-commuting variables, and in the last ten years, since the paper \cite{ando}, in the matrix analysis community, that is the current state of the art. The so called Karcher mean emerged among the different multivariable noncommutative geometric means as the most important one due to its connections to metric geometry as the center of mass, see for instance \cite{limpalfia}. Calculation of this particular mean is of interest in the numerical linear algebra community, see for instance \cite{Ba08,BMP,moakher} and many others. The best technique that we currently have relies on approximations using gradient descent and Newton-like iterations, see for example \cite{LP2} and the references therein.

On the other hand the theory of analytic functions in several non-commuting indeterminates were enjoying its own success, leading to recent representation theorems for free analytic functions \cite{helton0} with positive real part on the non-commutative operatorial unit $k$-ball, see for example \cite{popescu2,popescu,verbovetskyi} among the numerous contributions to this field. Another recent setting is considering free analytic functions mapping the upper operator poly-halfspace $\Pi^k$ into the upper operator halfspace $\Pi:=\{X\in\mathcal{B}(E):\Im{X}=\frac{1}{2i}(X-X^*)>0\}$ where $E$ is a Hilbert space and $k$ is a positive integer \cite{pascoe,agler0}. In this field operator monotonicity has been studied only recently in the case of commuting tuples of operators \cite{agler} and the focus remained restricted to free functions that are a priori assumed to be free analytic. In \cite{agler} an operator model theoretic representation formula is proved for $k$-variable continuously differentiable real functions that are operator monotone for commuting tuples of matrices. For a very recent contribution characterizing free analytic matrix monotone functions by establishing an operator model theoretic representation formula, we refer to \cite{pascoe} using techniques from \cite{agler} related to operator monotonicity and the so called Herglotz representation formula for free analytic functions with positive real part on the non-commutative operatorial $k$-ball \cite{popescu2,popescu}. The obtained representations in \cite{pascoe} based on \cite{popescu2} rely on the theory of noncommutative disc algebras developed for non-commutative power series representations and the Cuntz-Toeplitz algebra generated by partial isometries induced by the left- and right creation operators on the full Fock space. The free analytic functions are then represented by the slice map, or in other words, the conditional expectation of the inverse of certain linear combinations of tensor products on the Cuntz-Toeplitz algebra determined by a completely positive linear map which is itself uniquely determined by the values of the analytic function, for more technical details see \cite{popescu2,popescu}. There is also an earlier related result \cite{helton}, where it is proved that a free rational, hence automatically analytic, $k$-variable matrix convex function has a representation formula superficially similar to the ones in \cite{agler,pascoe}. The representation formula called the \emph{butterfly representation} provided in \cite{helton} for rational free analytic functions is the closest relative of our representation formula obtained in this paper.

Let $\mathbb{P}_n$ denote the cone of positive definite $n$-by-$n$ matrices, which is a subset of the cone of positive definite operators denoted by $\mathbb{P}(E)$ over the (infinite dimensional) Hilbert space $E$. In this paper we adopt the convention that if $X\in\mathbb{P}(E)$ then this means that $X$ is also lower bounded, hence invertible, the same applies to the notation $X>0$. We will denote by $\hat{\mathbb{P}}$ the cone of positive semi-definite operators and similarly $\hat{\mathbb{P}}_n$ denotes its finite $n$-by-$n$-dimensional counterpart. Let $\mathbb{S}(E)$ denote the vector space of bounded self-adjoint operators over $E$ and $\mathbb{S}_n$ its finite dimensional counterpart which is a subspace of $\mathbb{S}$. Sometimes we will denote $\mathbb{S}(E)$ simply by $\mathbb{S}$ if certain assertion holds for elements of $\mathbb{S}(E)$ for any  Hilbert space $E$. We will adopt the following shorthand notation for $k$-tuples of operators in the above sets:
\begin{equation*}
{X}:=(X_1,\ldots,X_k)
\end{equation*}
where $k$ is a fixed positive integer. We also use the notation ${I}:=(I,\ldots,I)$ where $I$ is the identity operator. The upper operator poly-halfspace is denoted by
\begin{equation*}
\Pi^k:=\{X\in\mathcal{B}(E)^k:\Im X_i>0, 1\leq i\leq k\}
\end{equation*}
where the \emph{imaginary part} of a bounded linear operator $A$ is defined as
\begin{equation*}
\Im A:=\frac{A-A^*}{2i}.
\end{equation*}
A several variable function $F:D\mapsto \mathbb{S}(E)$ for a domain $D\subseteq\mathbb{S}(E)^k$ defined for all Hilbert spaces $E$ is called a \textit{free} or \textit{noncommutative function} (NC function) if for all $A,B\in \mathbb{S}(E)^k$ which are in the domain of $F$ we have
\begin{itemize}
	\item[(1)] $F(U^{*}A_1U,\ldots,U^{*}A_kU)=U^{*}F(A_1,\ldots,A_k)U$ for all $U^{-1}=U^*\in\mathcal{B}(E)$,
	\item[(2)] $F\left(\left[ \begin{array}{cc}
A_1 & 0 \\
0 & B_1 \end{array} \right],\ldots,\left[ \begin{array}{cc}
A_k & 0 \\
0 & B_k \end{array} \right]\right)=\left[ \begin{array}{cc}
F(A_1,\ldots,A_k) & 0 \\
0 & F(B_1,\ldots,B_k) \end{array} \right]$,
\end{itemize}
Operator monotonicity and concavity is defined in Definition~\ref{D:opmon} and Definition~\ref{D:opconcave} below. It must be noted that the author does not know any operator mean in the literature which is operator monotone and does not satisfy the properties of free functions. In particular all the different geometric and operator means discussed in \cite{ando,bhatiaholbrook,BMP,lawsonlim0,lawsonlim,limpalfia,moakher,palfia,palfia2,palfia3} satisfy the free function property and of course monotonicity.

The main result of this paper, the representation formula equivalent to operator concavity, monotonicity and free analytic continuation to the upper operator poly-halfspace $\Pi^k$ mapping it to the upper operator halfspace $\Pi$ for a free function $F$ is as follows:
\begin{theoremintro}\label{T:0_Loewner_several_var}
Let $E$ be a Hilbert space and let $F:\mathbb{P}(E)^k\mapsto \mathbb{P}(E)$ be a free function. Then the following are equivalent:
\begin{itemize}
\item[(a)]$F$ is operator monotone;
\item[(b)]$F$ is operator concave;
\item[(c)]There exists a Hilbert space $\mathcal{K}$, a closed subspace $\mathcal{K}_0\leq \mathcal{K}$ and the corresponding orthogonal projection $P_{\mathcal{K}_0}$ with range $\mathcal{K}_0$, $B_i\in\hat{\mathbb{P}}(\mathcal{K})$, $0\leq i\leq k$ with $B_0\geq \sum_{i=1}^kB_i$ and a state $w\in\mathcal{B}_1^+(\mathcal{K})^*$ such that $F(X)$ is the $w$-conditional expectation of the Schur complement of the linear pencil $L_{B}(X)=B_0\otimes I+\sum_{i=1}^kB_i\otimes X_i$ pivoting on the subspace $\mathcal{K}_0$ for all $X\in\mathbb{P}(E)^k$. I.e. for all $X\in\mathbb{P}(E)^k$ we have
\begin{equation*}\label{eq:T:0_Loewner_several_var:0.1}
\begin{split}
F(X)=&w(B_{0,11})\otimes I+\sum_{i=1}^kw(B_{i,11})\otimes (X_i-I)\\
&-(w\otimes I)\left\{\left[B_{0,12}\otimes I+\sum_{i=1}^kB_{i,12}\otimes (X_i-I)\right]\right.\\
&\left[B_{0,22}\otimes I+\sum_{i=1}^kB_{i,22}\otimes (X_i-I)\right]^{-1}\\
&\left.\left[B_{0,21}\otimes I+\sum_{i=1}^kB_{i,21}\otimes (X_i-I)\right]\right\}
\end{split}
\end{equation*}
where
\begin{equation*}\label{eq:T:0_Loewner_several_var:0.2}
\begin{split}
{B}_{i,11}(A,v):=&P_{\mathcal{K}_0}{B}_{i}P_{\mathcal{K}_0},\\
{B}_{i,12}(A,v):=&P_{\mathcal{K}_0}{B}_{i}(I-P_{\mathcal{K}_0}),\\
{B}_{i,21}(A,v):=&(I-P_{\mathcal{K}_0}){B}_{i}P_{\mathcal{K}_0},\\
{B}_{i,22}(A,v):=&(I-P_{\mathcal{K}_0}){B}_{i}(I-P_{\mathcal{K}_0});
\end{split}
\end{equation*}
\item[(d)]$F$ has a free analytic continuation to $\Pi^k$, mapping $\Pi^k$ to $\Pi$.
\end{itemize}
\end{theoremintro}
In Theorem~\ref{T:0_Loewner_several_var} we restricted the domain and range of the free functions, but in principle M\"obius transformations can be used to transform the domain of operator monotone functions to $\mathbb{P}$ and the range assumption in Theorem~\ref{T:0_Loewner_several_var} is basically equivalent to the assumption that the function is bounded from below. The exact method to use M\"obius transformations on the domain of $F$ is discussed in the last section of the paper. Also due to the free analytic continuation part of Theorem~\ref{T:0_Loewner_several_var} the techniques in \cite{pascoe} apply, to obtain a general formula using the Cayley transform to connect $\Pi^k$ and the noncommutative operatorial $k$-ball and transform the noncommutative Herglotz formula established in \cite{popescu} for free holomorphic functions with real part on the noncommutative operatorial $k$-ball. This in principle follows the technique used in \cite{bhatia} to establish a Nevanlinna type formula for holomorphic functions mapping the upper complex half-plane into itself.

In order to prove the Theorem~\ref{T:0_Loewner_several_var}, we first establish the equivalence of operator monotonicity and operator concavity over certain domains. The proof of the equivalence between operator monotonicity and concavity in the one variable case has already appeared as early as in \cite{hansenpedersen,hansen3}. In particular the argument given in \cite{hansen3} goes through in several variables as well with minor modifications. Then by concavity we deduce the norm continuity of operator monotone functions. Then we introduce the hypograph of free functions as the saturation of their graphs:
$$\hypo(F)=(\hypo(F)(E)):=(\{(Y,X)\in\mathbb{S}(E)\times\mathbb{P}(E)^k:Y\leq F(X)\}).$$
Theorem~\ref{T:convex_hypographs} states that a free function $F$ is operator concave/operator monotone if and only if $\hypo(F)$ is a matrix convex set in the sense of Wittstock. Then in Proposition~\ref{P:separating_pencil} for each boundary point of $\hypo(F)(E)$ we construct a linear matrix pencil, similarly as in \cite{effros,helton3} for other matrix convex sets, that is positive semi-definite on the matrix convex set $\hypo(F)$ and it is singular at the boundary point. The proof of Proposition~\ref{P:separating_pencil} in principle is similar to the one given in section 6 of \cite{helton3}, but it also works for complex Hilbert spaces and it is formulated in such a way, that it can be applied directly in the infinite dimensional case as well, although for our purpose the finite dimensional version suffices.

The next step is Theorem~\ref{T:reconstruct_schur} which is an explicit extremal linear matrix inequality (LMI) solution formula based on the Schur complement that provides a key reconstruction formula for the actual values $F(X)v$ at arbitrary tuples $X$ applied to arbitrary unit vectors $v\in E$. Then we construct a non-separable Hilbert space where all possible combinations of the points of the domain of the free function $F$ are listed as a direct summand giving a block diagonal operator acting on this space. Then we apply Theorem~\ref{T:reconstruct_schur} to reconstruct the value of the function at this block diagonal operator acting on vectors of this non-separable Hilbert space to get Lemma~\ref{L:reconstruct_schur_multiple} and then apply a compactness argument to obtain the representation formula Theorem~\ref{T:Loewner_formula}. The compactness argument is inspired by the one given in \cite{helton} for proving linear dependence of free rational expressions over finite dimensional spaces. Theorem~\ref{T:Loewner_formula} expresses operator concave/operator monotone functions as a conditional expectation of the Schur complement of a linear matrix pencil establishing (c) in Theorem~\ref{T:0_Loewner_several_var}. Then in the last section of the paper we establish the free analytic continuation to obtain (d) in Theorem~\ref{T:0_Loewner_several_var}. Here we use Proposition~\ref{P:Drury}, a version of a result appeared in \cite{drury} for matrices, providing an estimate on the norm of the Schur complement of a sectorial operator.

The approach outlined above provides a new proof of L\"owner's theorem in the one variable case as well. The representation formula based on the Schur complement is also new even in this case. We believe that this approach in the one variable case has the advantage over the existing ones in that there is no need to establish continuous differentiability of operator monotone functions by mollifier smoothing techniques. Also the linear matrix pencil that appears in (c) in Theorem~\ref{T:0_Loewner_several_var} can be thought of as the direct sum of the LMI representations of the supporting linear functionals of the set $\hypo(F)$, so the formula in (c) is intuitive from the convex geometrical point of view.


\section{Noncommutative functions, monotonicity and concavity}

We put the following plausible assumptions on our operator valued several variable functions. Let $E$ denote an arbitrary Hilbert space.
\begin{definition}[NC function, \cite{verbovetskyi}]
A several variable function $F:D(E)\mapsto \mathbb{S}(E)$ for a domain $D(E)\subseteq\mathbb{S}(E)^k$ defined for all Hilbert spaces $E$ is called a \textit{free} or \textit{noncommutative function} (NC function) if for all $E$ and all $A,B\in D(E)\subseteq\mathbb{S}(E)^k$
\begin{itemize}
	\item[(1)] $F(U^{*}A_1U,\ldots,U^{*}A_kU)=U^{*}F(A_1,\ldots,A_k)U$ for all $U^{-1}=U^*\in\mathcal{B}(E)$,
	\item[(2)] $F\left(\left[ \begin{array}{cc}
A_1 & 0 \\
0 & B_1 \end{array} \right],\ldots,\left[ \begin{array}{cc}
A_k & 0 \\
0 & B_k \end{array} \right]\right)=\left[ \begin{array}{cc}
F(A_1,\ldots,A_k) & 0 \\
0 & F(B_1,\ldots,B_k) \end{array} \right]$.
\end{itemize}
Note that this already includes the closure of the domain $D(E)$ under direct sums and element-wise unitary conjugation. A free function $F$ can be regarded as a graded function $F:D(E)\mapsto \mathbb{S}(E)$ between the collection of domains $D=(D(E))$ and its range included in $\mathbb{S}(E)$ for each Hilbert space $E$.
\end{definition}

The second property means that NC functions respect direct sum decompositions, while the first property is invariancy under unitary conjugations. Both assumptions are plausible, since in one variables, the functional calculus has these properties. Also see \cite{agler,pascoe,popescu} which adopt the same assumptions, under the additional assumption of free analyticity. Moreover many such functions considered as means of operators have these properties \cite{ando,bhatia,bhatia2,bhatiaholbrook,BMP,lawsonlim12,lawsonleelim,limpalfia,moakher}.

The set $\mathbb{S}$ is equipped with a partial order, the positive definite order $\leq$, which means that for $A,B\in\mathbb{S}$ we have $A\leq B$ if and only if $0\leq B-A$, that is $\left\langle (B-A)x,x\right\rangle\geq 0$ for all vectors $x\in E$. For $k$-tuples $A,B\in\mathbb{S}^k$ we define similarly
\begin{equation}
A\leq B\text{ iff }A_i\leq B_i\text{ for all }i=1,\ldots,k.
\end{equation}

\begin{definition}[Monotonicity]\label{D:opmon}
An NC function $F:\mathbb{P}^k\mapsto \mathbb{S}$ is said to be operator monotone if whenever $A\leq B$ for $A,B\in\mathbb{P}^k$, we have
\begin{equation*}
F(A)\leq F(B).
\end{equation*}
If this property is satisfied only in finite dimensions, then we say that the NC function $F:\mathbb{P}_n^k\mapsto \mathbb{S}_n$ is $n$-monotone.
\end{definition}

\begin{definition}[Concavity \& Convexity]\label{D:opconcave}
A function $F:\mathbb{P}^k\mapsto \mathbb{S}$ is said to be operator concave if for all $A,B\in\mathbb{P}^k$ and $\lambda\in[0,1]$, we have
\begin{equation*}
(1-\lambda)F(A)+\lambda F(B)\leq F((1-\lambda)A+\lambda B)
\end{equation*}
Similarly if this property is satisfied only for $n$-by-$n$ matrices, then we say that the NC function $F:\mathbb{P}_n^k\mapsto \mathbb{S}_n$ is $n$-concave. Operator and $n$-convexity is defined accordingly. Many times later on we consider concavity or convexity on subsets of $\mathbb{P}^k$.
\end{definition}

All the above definitions can be considered on other convex domains (order intervals), not just $\mathbb{P}^k$.

\begin{remark}
Notice that if a function $F$ is $n$-monotone or $n$-concave, then by the direct sum property it is $m$-monotone and $m$-concave accordingly for all $1\leq m\leq n$. Also if $F$ is operator monotone or concave then it is $n$-monotone or $n$-concave for all finite $n\geq 1$. 
\end{remark}




We need the following version of Lemma 3.5.5. in \cite{niculescu}.
\begin{lemma}\label{L:concave-bounded}
Let $F$ be a concave function into $\mathbb{S}$ on an open convex set $U$ in a normed linear space. If $F$ is bounded from below in a neighborhood of one point of $U$, then $F$ is locally bounded on $U$.
\end{lemma}
\begin{proof}
Suppose that $F$ is bounded from below by $MI$ for some $M\in\mathbb{R}$ on an open ball $B(a,r)$ with radius $r$ around $a$. Let $x\in U$ and choose $\rho>1$ such that $z:=a+\rho(x-a)\in U$. If $\lambda=1/\rho$, then $V=\{v:v=(1-\lambda)y+\lambda z, y\in B(a,r)\}$ is a neighborhood of $x=(1-\lambda)a+\lambda z$, with radius $(1-\lambda)r$. Moreover, for $v\in V$ we have
\begin{equation*}
F(v)\geq (1-\lambda)F(y)+\lambda F(z)\geq (1-\lambda)MI+\lambda F(z)\geq KI
\end{equation*}
for some $K\in\mathbb{R}$. To show that $F$ is bounded above in the same neighborhood, choose arbitrarily $v\in V$ and notice that $2x-v\in V$. By concavity $F(x)\geq F(v)/2+F(2x-v)/2$, which yields
\begin{equation*}
F(v)\leq 2F(x)-F(2x-v)\leq 2F(x)-KI.
\end{equation*}
\end{proof}

We equip $\mathcal{B}(E)^k$ and similarly $\mathbb{S}^k$ and $\hat{\mathbb{P}}^k$ with the norm
\begin{equation*}
\|X\|:=\sum_{i=1}^k\|X_i\|
\end{equation*}
for a tuple $X\in\mathcal{B}(E)^k$.

\begin{proposition}[see also Proposition 3.5.4 in \cite{niculescu}]\label{P:concave-cont}
An operator concave function $F:\mathbb{P}^k\mapsto \mathbb{S}$ which is locally bounded from below, is continuous in the norm topology.
\end{proposition}
\begin{proof}
Let $U\subseteq\mathbb{P}^k$ be an open norm bounded neighborhood with respect to the operator norm $\|\cdot\|$. Let $A\in U$ and $r>0$ such that the open ball $B(A,2r):=\{X\in U:\|X-A\|<2r\}\subseteq U$. Let $X,Y\in B(A,r)$ and $X\neq Y$ such that $\alpha:=\|Y-X\|<r$. Let
\begin{equation}\label{eq:P:concave-cont-1}
Z:=Y+\frac{r}{\alpha}(Y-X).
\end{equation}
Then
\begin{equation*}
\|Z-A\|\leq\|Y-A\|+\frac{r}{\alpha}\|Y-X\|<2r,
\end{equation*}
i.e. $Z\in B(A,2r)$. By \eqref{eq:P:concave-cont-1} we have
\begin{equation*}
Y=\frac{r}{r+\alpha}X+\frac{\alpha}{r+\alpha}Z,
\end{equation*}
so by operator concavity of $F$ we get
\begin{equation*}
F(Y)\geq\frac{r}{r+\alpha}F(X)+\frac{\alpha}{r+\alpha}F(Z),
\end{equation*}
which after rearranging yields
\begin{equation*}
\begin{split}
F(X)-F(Y)&\leq\frac{\alpha}{r+\alpha}(F(X)-F(Z))\\
&\leq\frac{\alpha}{r+\alpha}2MI\leq\frac{\alpha}{r}2MI,
\end{split}
\end{equation*}
where the real number $M>0$ provides a local bound for $F$ on $U$ in the form of $-2MI\leq F(X)-F(Z)\leq 2MI$ by Lemma~\ref{L:concave-bounded}. Now exchange the role of $X$ and $Y$ in the above to obtain the reverse inequality
\begin{equation*}
F(Y)-F(X)\leq\frac{\alpha}{r}MI.
\end{equation*}
From the above pair of inequalities we get
\begin{equation*}
\|F(Y)-F(X)\|\leq\frac{M}{r}\|Y-X\|
\end{equation*}
proving continuity.
\end{proof}

A net of operators $\{A_i\}_{i\in \mathcal{I}}$ is called increasing if $A_i\geq A_j$ for $i\geq j$ and $i,j\in\mathcal{I}$. Also $\{A_i\}_{i\in \mathcal{I}}$ is bounded from above if there exists some real constant $K>0$ such that $A_i\leq KI$ for all $i\in\mathcal{I}$. It is well known that any bounded from above increasing net of operators $\{A_i\}_{i\in \mathcal{I}}$ has a least upper bound $\sup_{i\in\mathcal{I}}A_i$ such that $B_j:=A_j-\sup_{i\in\mathcal{I}}A_i$ converges to $0$ in the strong operator topology. Similarly if we have a decreasing net of bounded operators that is bounded from below, then the net converges to its greatest lower bound.

The next characterization result is an extension of Theorem 2.1 in \cite{hansen3} to several variables. The proof is analogous to that of Theorem 2.1. We consider the finite dimensional situation, but its proof is presented in such a way that it works also in the infinite dimensional setting as well.
\begin{proposition}\label{P:contmonotone}
Let $F:\mathbb{P}_{2n}^k\mapsto \mathbb{S}_{2n}$ be a $2n$-monotone function. Then its restriction $F:\mathbb{P}_{n}^k\mapsto \mathbb{S}_{n}$ is $n$-concave, moreover it is norm continuous.
\end{proposition}
\begin{proof}
Let $A,B\in\mathbb{P}_n^k$ and let $\lambda\in[0,1]$. Then the following $2n$-by-$2n$ block matrix is unitary
\begin{equation*}
V:=\left[ \begin{array}{cc}
\lambda^{1/2}I_n & -(1-\lambda)^{1/2}I_n \\
(1-\lambda)^{1/2}I_n & \lambda^{1/2}I_n \end{array} \right].
\end{equation*}
Elementary calculation reveals that
\begin{equation*}
V^*\left[ \begin{array}{cc}
A & 0 \\
0 & B \end{array} \right]V=
\left[ \begin{array}{cc}
\lambda A+(1-\lambda)B & \lambda^{1/2}(1-\lambda)^{1/2}(B-A) \\
\lambda^{1/2}(1-\lambda)^{1/2}(B-A) & (1-\lambda)A+\lambda B \end{array} \right].
\end{equation*}
Set $D:=-\lambda^{1/2}(1-\lambda)^{1/2}(B-A)$ and notice that for any given $\epsilon>0$
\begin{equation*}
\left[ \begin{array}{cc}
\lambda A+(1-\lambda)B+\epsilon I & 0 \\
0 & 2Z \end{array} \right]-V^*\left[ \begin{array}{cc}
A & 0 \\
0 & B \end{array} \right]V\geq
\left[ \begin{array}{cc}
\epsilon I & D \\
D & Z \end{array} \right]
\end{equation*}
if $Z \geq (1-\lambda)A+\lambda B$. The last $k$-tuple of block matrices is positive semi-definite if $Z \geq D_i^2/\epsilon$ for all $1\leq i\leq k$. So, for sufficiently large positive definite $Z$ we have
\begin{equation*}
V^*\left[ \begin{array}{cc}
A & 0 \\
0 & B \end{array} \right]V\leq \left[ \begin{array}{cc}
\lambda A+(1-\lambda)B+\epsilon I & 0 \\
0 & 2Z \end{array} \right].
\end{equation*}
For such $Z>0$, by the $2n$-monotonicity of $F$ we get
\begin{equation*}
F\left(V^*\left[ \begin{array}{cc}
A & 0 \\
0 & B \end{array} \right]V\right)\leq \left[ \begin{array}{cc}
F(\lambda A+(1-\lambda)B+\epsilon I) & 0 \\
0 & F(2Z) \end{array} \right].
\end{equation*}
We also have that
\begin{equation*}
\begin{split}
&F\left(V^*\left[ \begin{array}{cc}
{A} & 0 \\
0 & {B} \end{array} \right]V\right)=
V^*\left[ \begin{array}{cc}
F({A}) & 0 \\
0 & F({B}) \end{array} \right]V\\
&=\left[ \begin{array}{cc}
\lambda F({A})+(1-\lambda)F({B}) & \lambda^{1/2}(1-\lambda)^{1/2}(F({B})-F({A})) \\
\lambda^{1/2}(1-\lambda)^{1/2}(F({B})-F({A})) & (1-\lambda)F({A})+\lambda F({B}) \end{array} \right],
\end{split}
\end{equation*}
hence we obtain that
\begin{equation}\label{eq:P:concave}
\lambda F({A})+(1-\lambda)F({B})\leq F(\lambda{A}+(1-\lambda){B}+\epsilon I).
\end{equation}
Now since $F$ is $2n$-monotone, $F(X+\epsilon{I})$ for $\epsilon>0$ forms a decreasing net of operators bounded from below by $F(X)$, thus the right strong limit
$$F^{+}({X}):=\inf_{\epsilon>0}F(X+\epsilon{I})=\lim_{\epsilon\to 0+}F(X+\epsilon{I})$$
exists for all ${X}\in\mathbb{P}_{2n}^k$. Hence for any $\epsilon>0$, using \eqref{eq:P:concave}, we obtain
\begin{equation*}
\lambda F^{+}({A})+(1-\lambda)F^{+}({B})\leq \lambda F({A}+\epsilon{I})+(1-\lambda)F({B}+\epsilon{I})\leq F(\lambda{A}+(1-\lambda){B}+2\epsilon{I}).
\end{equation*}
Taking the strong limit $\epsilon\to 0+$ we obtain that
\begin{equation*}
\lambda F^{+}({A})+(1-\lambda)F^{+}({B})\leq F^{+}(\lambda{A}+(1-\lambda){B}),
\end{equation*}
i.e. the NC function $F^{+}$ is $n$-concave. Also
\begin{equation*}
F({X})\leq F^{+}({X})\leq F({X}+\epsilon{I})
\end{equation*}
for all $\epsilon>0$, since $F$ is monotone increasing, hence $F^{+}$ is bounded from below on order bounded sets, so by Proposition~\ref{P:concave-cont} $F^{+}$ is norm continuous on order bounded sets, since every point $A\in\mathbb{S}$ has a basis of neighborhoods in the norm topology that are order bounded sets.
As the last step, again by the monotonicity of $F$ we have
\begin{equation*}
F^{+}({X}-\epsilon{I})\leq F({X})\leq F^{+}({X}),
\end{equation*}
and since $F^{+}$ is norm continuous we get $F=F^{+}$ by taking the norm limit $\epsilon\to 0+$. Hence we can also take the norm limit $\epsilon\to 0+$ in \eqref{eq:P:concave} and conclude that $F$ is $n$-concave and continuous.
\end{proof}

Since the above proof also works in infinite dimensions we have the following result.
\begin{corollary}\label{P:concavemonotone}
An operator monotone NC function $F:\mathbb{P}^k\mapsto \mathbb{S}$ is operator concave and norm continuous.
\end{corollary}

The reverse implication is also true if $F$ is bounded from below:
\begin{theorem}\label{T:concavemonotone}
Let $F:\mathbb{P}^k\mapsto \hat{\mathbb{P}}$ be operator concave ($n$-concave) NC function. Then $F$ is operator monotone ($n$-monotone).
\end{theorem}
\begin{proof}
The proof goes along the lines of Theorem 2.3 in \cite{hansen3}.
\end{proof}

By the above if we wish to characterize operator monotone functions $F:\mathbb{P}^k\mapsto \mathbb{P}$, then it suffices to characterize operator concave ones.

We close the section with a semi-continuity property for operator monotone NC functions.

\begin{proposition}\label{P:upper_semi_continuity}
Let $F:\mathbb{P}^k\mapsto \hat{\mathbb{P}}$ be an operator monotone NC function. Then for any bounded from above increasing net of $k$-tuple of operators $\{A_i\}_{i\in \mathcal{I}}$ with $A_i\in\mathbb{P}^k$ we have that $F$ is strongly upper semi-continuous along increasing nets, i.e.
\begin{equation*}
\sup_{i\in\mathcal{I}}F(A_i)\leq F\left(\sup_{i\in\mathcal{I}}A_i\right).
\end{equation*}
\end{proposition}
\begin{proof}
We have that $A_j\leq \sup_{i\in\mathcal{I}}A_i$ component-wise for all $j\in\mathcal{I}$, i.e. $(A_j)_l\leq \sup_{i\in\mathcal{I}}(A_i)_l$ for all $1\leq l\leq k$ and $j\in\mathcal{I}$. Thus by the monotonicity of $F$ we have
\begin{equation*}
F(A_j)\leq F\left(\sup_{i\in\mathcal{I}}A_i\right)
\end{equation*}
for all $j\in\mathcal{I}$, so it follows that
\begin{equation*}
\sup_{j\in\mathcal{I}}F(A_j)\leq F\left(\sup_{i\in\mathcal{I}}A_i\right).
\end{equation*}
\end{proof}

\section{Supporting linear pencils and hypographs}
In this section we will use the theory of matrix convex sets introduced first by Wittstock. Some references on free convexity and matrix convex sets are \cite{effros,helton,helton2,helton3,helton4}. Let $\Lat(E)$ denote the \emph{lattice of subspaces} of $E$. The notation $K\leq E$ means that $K$ is a closed subspace of $E$, hence a Hilbert space itself.
\begin{definition}[Matrix/Freely convex set]\label{def:matrix_covexity}
A graded collection $C=(C(K))$, where each $C(K)\subseteq \mathbb{S}(K)^k$ and $K$ is a Hilbert space, is a bounded $\tau$ open/closed \emph{matrix convex} or \emph{freely convex} set if
\begin{itemize}
\item[(i)] each $C(K)$ is open/closed in the $\tau$ topology;
\item[(ii)] $C$ respects direct sums, i.e. if $(X_1,\ldots,X_k)\in C(K)$ and $(Y_1,\ldots,Y_k)\in C(N)$ and $Z_j:=\left[ \begin{array}{cc}
X_j & 0 \\
0 & Y_j \end{array} \right]$ for Hilbert spaces $K,N$, then $(Z_1,\ldots,Z_k)\in C(K\oplus N)$;
\item[(iii)] $C$ respects conjugation with isometries, i.e. if $Y\in C(K)$ and $V:N\mapsto K$ is an isometry for Hilbert spaces $K,N$, then $V^*YV=(V^*Y_1V,\ldots,V^*Y_kV)\in C(N)$;
\item[(iv)] each $C(K)$ is bounded.
\end{itemize}
\end{definition}

The above definition has some equivalent characterizations under slight additional assumptions.

\begin{definition}
A graded collection $C=(C(K))$, where each $C(K)\subseteq \mathbb{S}(K)^k$, is \emph{closed with respect to reducing subspaces} if for any tuple of operators $(X_1,\ldots,X_k)\in C(K)$ and any corresponding mutually invariant subspace $N\subseteq K$, we have that the restricted tuple $(\hat{X}_1,\ldots,\hat{X}_k)\in C(N)$, where each $\hat{X}_i$ is the restriction of $X_i$ to the invariant subspace $N$ for all $1\leq i\leq k$.
\end{definition}

\begin{lemma}[Lemma 2.3 in \cite{helton4}, \S 2 in \cite{helton2}]\label{lem:matrix_convexity_equivalent}
Suppose that the graded collection $C=(C(K))$, where each $C(K)\subseteq \mathbb{S}(K)^k$ respects direct sums in the sense as in (ii) in Definition~\ref{def:matrix_covexity} and it respects unitary conjugation in the sense as in (iii) in Definition~\ref{def:matrix_covexity} with $N=K$.
\begin{itemize}
\item[(1)] If $C$ is closed with respect to reducing subspaces then $C$ is matrix convex if and only if each $C(K)$ is convex in the classical sense of taking scalar convex combinations.
\item[(2)] If $C$ is (nonempty and) matrix convex, then  $0=(0,\ldots,0)\in C(1)$ if and only if $C$ is closed with respect to simultaneous conjugation by contractions.
\end{itemize}
\end{lemma}


Given a set $A\subseteq\mathbb{S}$ we define its \emph{saturation} as $$\sat(A):=\{X\in\mathbb{S}:\exists Y\in A, Y\geq X\}.$$
Similarly for a graded collection $C=(C(K))$, where each $C(K)\subseteq\mathbb{S}(K)$, its \emph{saturation} $\sat(C)$ is the disjoint union of $\sat(C(K))$ for each Hilbert space $K$.

\begin{definition}[Hypographs]
Let $F:\mathbb{P}^k\mapsto \mathbb{S}$ be an NC function. Then for a fixed constant $c>0$, we define its \emph{hypograph} $\hypo(F)$ as the graded collection of the saturation of its image, i.e.
$$\hypo(F)=(\hypo(F)(K)):=(\{(Y,X)\in\mathbb{S}(K)\times\mathbb{P}(K)^k:Y\leq F(X)\}),$$
\end{definition}

\begin{theorem}\label{T:convex_hypographs}
Let $F:\mathbb{P}^k\mapsto \mathbb{S}$ be an NC function. Then its hypograph $\hypo(F)$ is a matrix convex set if and only if $F$ is operator concave.
\end{theorem}
\begin{proof}
Suppose first that $F$ is operator concave. We will prove the matrix convexity of $\hypo(F)$ by establishing the properties in (1) of Lemma~\ref{lem:matrix_convexity_equivalent}. By the definition of operator concavity and the convexity of $\mathbb{P}$ and the order intervals, it follows easily that for each Hilbert space $K$, $\hypo(F)(K)$ is convex in the usual sense of taking scalar convex combinations. To see that $\hypo(F)$ is closed with respect to reducing subspaces, assume that $(Y,X)\in\hypo(F)(L)$ with $(Y,X)=(\hat{Y},\hat{X})\oplus(\overline{Y},\overline{X})$ and $(\hat{Y},\hat{X})\in\mathbb{S}(K)\times\mathbb{P}(K)^k, (\overline{Y},\overline{X})\in\mathbb{S}(N)\times\mathbb{P}(N)^k$ for Hilbert spaces $K\oplus N=L$. Then since $F$ is an NC function, it respects direct sums, hence $Y\leq F(X)=F(\hat{X})\oplus F(\overline{X})$. Again by the definition of NC functions, we have $F(\hat{X})\in\mathbb{S}(K)$ and $F(\overline{X})\in\mathbb{S}(N)$. Since $Y=\hat{Y}\oplus\overline{Y}$, it follows that $\hat{Y}\leq F(\hat{X})$ and $\overline{Y}\leq F(\overline{X})$, i.e. $(\hat{Y},\hat{X})\in\hypo(F)(K)$ and $(\overline{Y},\overline{X})\in\hypo(F)(N)$.

For the converse, suppose that $\hypo(F)$ is a matrix convex set. First of all notice that $\hypo(F)$ is closed with respect to reducing subspaces. Indeed, similarly to the above assume that $(Y,X)\in\hypo(F)(L)$ with $(Y,X)=(\hat{Y},\hat{X})\oplus(\overline{Y},\overline{X})$ and $(\hat{Y},\hat{X})\in\mathbb{S}(K)\times\mathbb{P}(K)^k, (\overline{Y},\overline{X})\in\mathbb{S}(N)\times\mathbb{P}(N)^k$ for Hilbert spaces $K\oplus N=L$. Then since $F$ is an NC function, it respects direct sums, hence $Y\leq F(X)=F(\hat{X})\oplus F(\overline{X})$. Again by the definition of NC functions, we have $F(\hat{X})\in\mathbb{S}(K)$ and $F(\overline{X})\in\mathbb{S}(N)$. Since $Y=\hat{Y}\oplus\overline{Y}$, it follows that $\hat{Y}\leq F(\hat{X})$ and $\overline{Y}\leq F(\overline{X})$, i.e. $(\hat{Y},\hat{X})\in\hypo(F)(K)$ and $(\overline{Y},\overline{X})\in\hypo(F)(N)$. So, again by (1) of Lemma~\ref{lem:matrix_convexity_equivalent} it follows that for each Hilbert space $L$, $\hypo(F)(L)$ is convex in the usual sense. Notice also that for each $L$ we can recover the values of the NC function $F$, since for a fixed $X\in\mathbb{P}(L)^k$ we have that $F(X)=\sup\{Y\in\mathbb{S}(L):(Y,X)\in\hypo(F)(L)\}$. In other words for all $t\in[0,1]$ and $A,B\in\mathbb{P}(L)^k$ we have that the tuple $(Y,X):=(1-t)(F(A),A)+t(F(B),B)$ is in $\hypo(F)(L)$, moreover since $F(X)=\sup\{Y\in\mathbb{S}(L):(Y,X)\in\hypo(F)(L)\}$ we have that $(1-t)F(A)+tF(B)\leq F(X)=F((1-t)A+tB)$ for all $A,B\in\mathbb{P}(L)^k$ and Hilbert space $L$, hence $F$ is operator concave.
\end{proof}

The above Theorem~\ref{T:convex_hypographs} combined with Theorem~\ref{T:concavemonotone} leads to the following:

\begin{corollary}\label{C:convex_hypographs_monotone}
Let $F:\mathbb{P}^k\mapsto \hat{\mathbb{P}}$ be an NC function. Then its hypograph $\hypo(F)$ is a matrix convex set if and only if $F$ is operator monotone.
\end{corollary}

A sharpening of Theorem~\ref{T:convex_hypographs} is possible if we establish further continuity properties of NC functions. We will use the terminology of \emph{measurable domain}, \emph{compact domain} and \emph{measurable operator function}, \emph{continuous operator function} given in \cite{kruszynski} in the proof of the following auxiliary result.

\begin{lemma}\label{L:strong_continuity}
Suppose that $F:\mathbb{P}(E)^k\mapsto \mathbb{S}(E)$ is a norm continuous NC function for any separable Hilbert space $E$. Then $F$ is continuous in the strong operator topology on bounded sets $\{X\in\mathbb{P}(E)^k:cI\leq X_i\leq CI\}$ for fixed constants $0<c<C$.
\end{lemma}
\begin{proof}
The set $D:=\{X\in\mathbb{P}(E)^k:cI\leq X_i\leq CI\}$ for fixed constants $0<c<C$ and any separable Hilbert space $E$ is a measurable domain. Moreover since $D$ is closed in the norm topology for fixed $E$ by Theorem 1.3. \cite{kruszynski} $D$ is a compact domain. Since $F$ is a measurable norm continuous operator function defined on the compact domain $D$ it follows by Theorem 1.6. 2) \cite{kruszynski} that $F$ is a continuous operator function. By the non-commutative version of the Stone-Weierstrass Theorem 3.4. \cite{kruszynski} the $C^*$-algebra of continuous operator functions on a compact domain are generated by the coordinate functions, i.e. the non-commutative polynomials are $C^*$-norm dense in the $C^*$-algebra of continuous operator functions. Hence it follows that $F$ can be approximated uniformly by non-commutative polynomials in the norm topology. Since operator multiplication is jointly strong operator continuous on norm bounded sets it follows that any non-commutative polynomial is jointly strong operator continuous on norm bounded sets. Now since $F$ can be approximated uniformly by strong operator continuous non-commutative polynomials in the norm topology, hence uniformly in the strong operator topology on $D$, it follows that $F$ is strong operator continuous on $D$.
\end{proof}

\begin{remark}
In \cite{kruszynski} the above used theorems are proved for separable Hilbert spaces, since the constructions use sets like the set of all Hilbert spaces, which does not exist, if we also include non-separable spaces. Hence one must restrict to separable Hilbert spaces (which as a matter of fact are isomorphic to each other for a fixed dimension) to obtain the results in \cite{kruszynski}.
\end{remark}

Given a Hilbert space $E$, the space of bounded linear operators $\mathcal{B}(E)$ as a von Neumann algebra has a unique predual $\mathcal{B}(E)_{*}$, the Banach space of trace-class operators. The topology induced by the duality $(\mathcal{B}(E),\mathcal{B}(E)_{*})$ is the $\sigma$- or \emph{ultra-weak} operator topology. In other words this topology is generated by the closed subspace $\mathcal{B}(E)_{*}$ of normal linear functionals of the dual space $\mathcal{B}(E)^{*}$.

\begin{theorem}\label{T:weakly_convex_hypographs}
Let $E$ be separable and let $F:\mathbb{P}^k\mapsto \mathbb{P}$ be an NC function. Then its restricted hypograph
\begin{equation}
\hypo_{c,C}(F)(E):=(\{(Y,X)\in\mathbb{S}(E)\times\{X\in\mathbb{P}(E)^k:cI\leq X_i\leq CI\}:Y\leq F(X)\})
\end{equation}
for any given $C>c>0$, is a $\sigma$-strongly/$\sigma$-weakly closed matrix convex set if and only if $F$ is operator concave if and only if $F$ is operator monotone.
\end{theorem}
\begin{proof}
By Theorem~\ref{T:convex_hypographs} and Corollary~\ref{C:convex_hypographs_monotone} $F$ is operator concave if and only if $F$ is operator monotone, if and only if $\hypo(F)(E)$ is a matrix convex set. Then also by Corollary~\ref{P:concavemonotone} $F$ is norm continuous, hence by Lemma~\ref{L:strong_continuity} $F$ is strong operator continuous on norm bounded sets, and also $\{X\in\mathbb{P}(E)^k:cI\leq X_i\leq CI\}$ is a strong operator closed convex set. Hence $\hypo_{c,C}(F)(E)$ is also strong operator closed and convex for each separable $E$. Since $\hypo_{c,C}(F)(E)$ is a strong operator closed convex set, its closure in the $\sigma$-weak operator topology is itself, since the weak and strong operator closure of convex sets of operators are the same and the weak and $\sigma$-weak operator topologies coincide on norm bounded sets, see Theorem II.2.6 in \cite{takesaki}.
\end{proof}

In the case of non-separable $E$, we cannot use Lemma~\ref{L:strong_continuity} to prove strong continuity of norm continuous functions on norm bounded sets. Instead we will consider some additional plausible assumption on $F$ in the case of non-separable $E$.

\begin{assumption}\label{A:semi-continuity2}
Let $F:\mathbb{P}^k\mapsto \mathbb{S}$ be an operator monotone NC function. If $E$ is non-separable we assume that for any bounded from above increasing net of $k$-tuple of operators $\{A_i\}_{i\in \mathcal{I}}$ with $A_i\in\mathbb{P}(E)^k$ we have that
\begin{equation}\label{eq:A:semi-continuity:2}
\sup_{i\in\mathcal{I}}F(A_i)\geq F\left(\sup_{i\in\mathcal{I}}A_i\right).
\end{equation}
\end{assumption}
For example the Karcher mean of positive definite operators satisfies the above assumptions, see \cite{lawsonlim1}, as well as any Kubo-Ando mean \cite{kubo} or the operator means in \cite{palfia2}.

In what follows we will consider supporting linear pencils for a matrix convex set, that are in one to one correspondence with supporting linear functionals coming from the Hahn-Banach Theorem given for topological vector spaces.

\begin{definition}[free $\epsilon$-neighborhood]
Given $\epsilon>0$ the \emph{free $\epsilon$-neighborhood of $0$}, denoted by $\mathfrak{N}_\epsilon$, is for each Hilbert space $K$, the graded collection $(\mathfrak{N}_\epsilon(K))$ where
$$\mathfrak{N}_\epsilon(K):=\{X\in\mathbb{S}(K)^k:\sum^k_{j=1}\|X_j\|<\epsilon\}.$$
\end{definition}

\begin{definition}[linear pencil]
A \emph{linear pencil} is an expression of the form
$$L(x):=A_0+A_1x_1+\cdots+A_kx_k$$
where each $A_i\in\mathbb{S}(K)$ for some Hilbert space $K$ whose dimension is the \emph{size} of the pencil $L$. The pencil is \emph{monic} if $A_0=I$ and then $L$ is a \emph{monic linear pencil}. We extend the evaluation of $L$ from scalars to operators by tensor multiplication. In particular $L$ evaluates at a tuple $X\in\mathbb{S}(N)^k$ as
$$L(X):=A_0\otimes I+A_1\otimes X_1+\cdots+A_k\otimes X_k.$$
We then regard $L(X)$ as a self-adjoint element of $\mathbb{S}(K\otimes N)$.
\end{definition}

Let $\mathcal{B}^+_1(K)_{*}\subset \mathcal{B}_1(K)_{*}$ denote the convex set of positive semi-definite operators over $K$ of trace one and let $\mathcal{B}^+_1(K)^{*}$ denote the state space of the $C^*$-algebra $\mathcal{B}(K)$. Note that positive linear functionals on unital $C^*$-algebras attain their norm at the unit, hence $\mathcal{B}^+_1(K)^{*}$ is convex and weak-$*$ compact by Banach-Alaoglu. Each element $T\in\mathcal{B}^+_1(K)_{*}$ corresponds to a state on $\mathbb{S}(K)$ by
$$X\mapsto \mathrm{tr}(XT)$$
for $X\in\mathbb{S}(K)$. We will need a version of Proposition 6.4 in \cite{helton3}, before giving our version we state the following auxiliary result from \cite{helton3}. We endow $\mathcal{B}^+_1(K)^{*}$ with the relative weak $*$-topology induced by the duality $(\mathcal{B}^h(K),\mathcal{B}^h(K)^{*})$, where the superscript $h$ denotes the self-adjoint part. 

\begin{lemma}\label{L:existTop}
Suppose $\mathcal{F}$ is a convex set of weak-$*$ continuous affine linear mappings $f:\mathcal{B}^+_1(K)^{*}\mapsto \mathbb{R}$ with respect to a duality. If for each $f\in\mathcal{F}$ there is a $T\in\mathcal{B}^+_1(K)^{*}$ such that $f(T)\geq 0$, then there is a $\mathcal{T}\in\mathcal{B}^+_1(K)^{*}$ such that $f(\mathcal{T})\geq 0$ for every $f\in\mathcal{F}$.
\end{lemma}
\begin{proof}
For $f\in\mathcal{F}$, let
\begin{equation*}
B_f:=\{T\in\mathcal{B}^+_1(K)^{*}:f(T)\geq 0\}\subset\mathcal{B}^+_1(K)^{*}.
\end{equation*}
By hypothesis, each $B_f$ is nonempty and it suffices to prove that
\begin{equation*}
\cap_{f\in\mathcal{F}}B_f\neq\emptyset.
\end{equation*}
Since each $B_f$ is compact, it suffices to prove that the collection $\{B_f:f\in\mathcal{F}\}$ has the finite intersection property. Let $f_1,\ldots,f_m\in\mathcal{F}$ be given. Suppose that
\begin{equation}\label{eq:L:existTop:1}
\cap_{j=1}^mB_{f_n}=\emptyset.
\end{equation}
Define $F:\mathcal{B}^+_1(K)^{*}\mapsto\mathbb{R}^m$ by
\begin{equation*}
F(T):=(f_1(T),\ldots,f_m(T)).
\end{equation*}
Then $F(\mathcal{B}^+_1(K)^{*})$ is both convex and compact because $\mathcal{B}^+_1(K)^{*}$ is and each $f_j$, hence $F$, is weak-$*$ continuous affine linear. Moreover $F(\mathcal{B}^+_1(K)^{*})$ does not intersect
\begin{equation*}
\mathbb{R}^m_+=\{x=(x_1,\ldots,x_m):x_j\geq 0\text{ for each }j\}.
\end{equation*}
Hence by the Hahn-Banach theorem there exists a linear functional $\lambda:\mathbb{R}\mapsto\mathbb{R}$ such that $\lambda(F(\mathcal{B}^+_1(K)^{*}))<0$ and $\lambda(\mathbb{R}^m_+)\geq 0$. We can write $\lambda$ as $\lambda(x)=\sum_{i=1}^m\lambda_jx_j$. Since $\lambda(\mathbb{R}^m_+)\geq 0$, it follows that each $\lambda_j\geq 0$. We have $\lambda_j\neq 0$ for at least one $1\leq j\leq m$, so without loss of generality we can assume that $\sum_{j=1}^m\lambda_j=1$. Let
\begin{equation*}
f:=\sum_{j=1}^m\lambda_jf_j.
\end{equation*}
Since $\mathcal{F}$ is convex, we have $f\in\mathcal{F}$. On the other hand, $f(T)=\lambda(F(T))$, hence if $T\in\mathcal{B}^+_1(K)^{*}$, then $f(T)<0$. Thus, for this $f$ there does not exist a $T\in\mathcal{B}^+_1(K)^{*}$ such that $F(T)\geq 0$, contradicting \eqref{eq:L:existTop:1}.
\end{proof}

\begin{lemma}\label{L:existT}
Let $C=(C(K))$ be a matrix convex set, where $C(K)\subseteq \mathbb{S}(K)^k$ and $(0,\ldots,0)\in C(\mathbb{C})$. Let a linear functional $\Lambda:\mathbb{S}(N)^k\mapsto\mathbb{R}$ be given for a fixed $N$ with $\dim(N)<\infty$. If $\Lambda(X)\leq 1$ for each $X\in C(N)$, then there exists a $T\in\mathcal{B}^+_1(N)^{*}$ such that for each Hilbert space $K$, and each $Y\in C(K)$ and each $V:N\mapsto K$ contraction $V$ we have
$$\Lambda(V^*YV)\leq \tr(VTV^*).$$
\end{lemma}
\begin{proof}
Since $\dim(N)<\infty$ it follows that $\mathcal{B}^+_1(N)_{*}=\mathcal{B}^+_1(N)^{*}$. For a Hilbert space $K$, a tuple $Y\in C(K)$ and a $V:N\mapsto K$ contraction, define $f_{Y,V}:\mathcal{B}^+_1(N)^{*}\mapsto\mathbb{R}$ by
\begin{equation*}
f_{Y,V}(T):=\tr(VTV^*)-\Lambda(V^*YV).
\end{equation*}
We claim that the collection $\mathcal{F}:=\{f_{Y,V}:Y,V\}$ is a convex set. Let $\lambda_i\geq 0$ for $1\leq i\leq n$ for a fixed integer $n$ and let $\sum_{i=1}^n\lambda_i=1$. Also let $(Y_i,V_i)$ be given where $Y_i\in C(K_i)$ for a Hilbert space $K_i$ and $V_i:N\mapsto K_i$ be a contraction for each $1\leq i\leq n$. Let $Z:=\oplus_{i=1}^nY_i$ and let $F$ denote the column operator matrix with entries $\sqrt{\lambda_i}V_i$. Then $Z\in C(\oplus K_i)$ and
\begin{equation*}
F^*F=\sum_{i=1}^n\lambda_iV_i^*V_i\leq \sum_{i=1}^n\lambda_iI=I.
\end{equation*}
By definition
\begin{equation*}
\sum_{i=1}^n\lambda_iV_i^*Y_iV_i=F^*ZF
\end{equation*}
and
\begin{equation*}
\sum_{i=1}^n\lambda_i\tr(V_iTV_i^*)=\tr(FTF^*)
\end{equation*}
for $T\in\mathcal{B}^+_1(N)^{*}$. Hence
\begin{equation*}
\sum_{i=1}^n\lambda_if_{Y_i,V_i}(T)=f_{Z,F}(T).
\end{equation*}

If $V$ has operator norm 1, we can choose a pure state $\gamma^*(\cdot)\gamma$ where $\gamma$ is a unit vector in $N$ such that
\begin{equation*}
1=\|V\gamma\|^2=\gamma^*V^*V\gamma=\tr(\gamma\gamma^*V^*V)=\tr(V\gamma\gamma^*V^*).
\end{equation*}
Then for $T=\gamma\gamma^*$ it follows that
\begin{equation*}
f_{Y,V}(T)=\tr(VTV^*)-\Lambda(V^*YV)=1-\Lambda(V^*YV).
\end{equation*}
Since $V^*YV\in C(N)$, the right hand side above is nonnegative. If the contraction $V$ does not have norm one, we can rescale it to have norm 1 and follow the same argument to show that $f_{Y,V}(T)\geq 0$. So, for each $f_{Y,V}$ there exists a $T\in\mathcal{B}^+_1(N)^{*}$ such that $f_{Y,V}(T)\geq 0$, moreover each $f_{Y,V}$ is weak-$*$ continuous. From Lemma~\ref{L:existTop} there exists a $\mathcal{T}\in\mathcal{B}^+_1(N)^{*}$ such that $f_{Y,V}(\mathcal{T})\geq 0$ for every $Y$ and $V$.
\end{proof}

For an arbitrary set $S$ of a topological vector space we will denote by $S^\circ$ its interior which is the union of all open sets contained in $S$. Below is a modified and generalized version of Proposition 6.4 in \cite{helton3}. A similar result is Theorem 5.4 in \cite{effros}.



\begin{proposition}\label{P:separating_pencil}
Let $F:\mathbb{P}^k\mapsto \mathbb{P}$ be an operator monotone function and let $N$ be a fixed Hilbert space with $\dim(N)<\infty$.
Then for each $A\in\mathbb{P}(N)^k$ and each unit vector $v\in N$ there exists a linear pencil
\begin{equation*}
L_{F,A,v}(Y,X):=B(F,A,v)_0\otimes I-vv^*\otimes Y+\sum_{i=1}^kB(F,A,v)_i\otimes (X_i-I)
\end{equation*}
of size $\dim(N)$ which satisfies the following properties:
\begin{itemize}
\item[(a)]$B(F,A,v)_i\in\mathcal{B}^+(N)_{*}$ and $\sum_{i=1}^kB(F,A,v)_i\leq B(F,A,v)_0$;
\item[(b)]For all $(Y,X)\in\hypo(F)$ we have $L_{F,A,v}(Y,X)\geq 0$;
\item[(c)]If $c_1I\leq A_i\leq c_2I$ for all $1\leq i\leq k$ and some fixed real constants $c_2>c_1>0$, then $\tr\{B(F,A,v)_0\}\leq \frac{F(c_2,\ldots,c_2)}{\min(1,c_1)}$.
\end{itemize}
\end{proposition}
\begin{proof}
By Theorem~\ref{T:convex_hypographs} $\hypo(F)$ is a matrix convex set. Consider the translated set $H(K):=\{(Y,X)\in\mathbb{S}(K)\times\mathbb{S}(K)^k:(Y,(X_1+I,\ldots,X_k+I))\in\hypo(F)(K)\}$. Still $H=(H(K))$ is a matrix convex set. Moreover since $F$ is positive we have $F(s I,\ldots,s I)=cI$ for some arbitrary small, but fixed $s>0$ and for any $X\geq (s I,\ldots,s I)$ we have $F(X)\geq cI$ by operator monotonicity. Hence $\hypo(F)$ contains a free $\epsilon$-neighborhood of $1$, so $H$ 
 contains a free $\epsilon$-neighborhood of $0$ for small enough $\epsilon>0$. Then $(F(A),A)$ is in the boundary of $\hypo(F)(N)$, hence $(F(A),(A_1-I,\ldots,A_k-I))$ is in the boundary of $H(N)$. Consider the real valued function $h:[\mathbb{P}(N)-I]^k\mapsto\mathbb{R}$ defined by $h(X):=v^*F(X+I)v$. Since $F$ is an operator concave function it follows that $h$ is concave and by Corollary~\ref{P:concavemonotone} it is also continuous in the norm topology.
It follows from the supporting hyperplane version of the Hahn-Banach theorem for the finite dimensional vector space $\mathbb{R}\times\mathbb{S}(N)^k$, more precisely Theorem 7.12 and 7.16 \cite{aliprantis} that the norm continuous convex function $g(X):=-h(X)=-v^*F(X+I)v$ has a subgradient at each interior point of its domain, hence at $(A-I)$ for $A\in\mathbb{P}(N)^k$. I.e. there exists a norm continuous linear functional $\lambda$ in the dual space of $\mathbb{S}(N)^k$ such that
\begin{equation*}
g(X)\geq g(A-I)+\lambda(X-A+I)
\end{equation*}
for all $X\in[\mathbb{P}(N)-I]^k$. Hence it follows that there exists $c\in\mathbb{R}$ and $l_i$ in the dual space of $\mathbb{S}(N)$ equipped with the norm topology such that
\begin{equation}\label{eq:P:separating_pencil:1}
1\geq \frac{1}{c}v^*F(X)v-\sum_{i=1}^kl_i(X_i-I)
\end{equation}
for all $X\in\mathbb{P}(N)^k$ and
\begin{equation}\label{eq:P:separating_pencil:1.2}
1=\frac{1}{c}v^*F(A)v-\sum_{i=1}^kl_i(A_i-I).
\end{equation}
Without loss of generality we can assume that $c>0$, hence it follows from \eqref{eq:P:separating_pencil:1} and the definition of $\hypo(F)(N)$ that
\begin{equation}\label{eq:P:separating_pencil:2}
1\geq \frac{1}{c}v^*Yv-\sum_{i=1}^kl_i(X_i-I)
\end{equation}
for all $(Y,X)\in\hypo(F)(N)$. From Lemma~\ref{L:existT} and \eqref{eq:P:separating_pencil:2} there exists a $T\in\mathcal{B}^+_1(N)_{*}$ such that for each Hilbert space $K$, and each $(Y,X)\in \hypo(F)(K)$ and each $V:N\mapsto K$ contraction $V$ we have
\begin{equation}\label{eq:P:separating_pencil:3}
\tr(VTV^*)-\frac{1}{c}v^*V^*YVv+\sum_{i=1}^kl_i(V^*(X_i-I)V)\geq 0.
\end{equation}
Since the dual space of $\mathbb{S}(N)$ (equipped with the norm topology) is the space of self-adjoint trace-class operators over $N$ we have $l_i(Z)=\tr\{B_iZ\}$ for all $Z\in\mathbb{S}(N)$ where $B_i\in\mathcal{B}^h(N)_{*}$ is a trace class operator. Hence we can write \eqref{eq:P:separating_pencil:3} as
\begin{equation}\label{eq:P:separating_pencil:4}
\tr(VTV^*)-\frac{1}{c}\tr\{vv^*V^*YV\}+\sum_{i=1}^k\tr\{B_i(V^*(X_i-I)V)\}\geq 0,
\end{equation}
moreover \eqref{eq:P:separating_pencil:1.2} becomes
\begin{equation}\label{eq:P:separating_pencil:4.1}
\tr(T)+\sum_{i=1}^k\tr\{B_i(A_i-I)\}=\frac{1}{c}\tr\{vv^*F(A)\}.
\end{equation}

Let $L_{B,v}$ denote the linear pencil
$$L_{B,v}(Y,X):=vv^*\otimes Y-\sum_{i=1}^kcB_i\otimes (X_i-I).$$
Let $K$ be a Hilbert space, let $\{e_i\}_{i\in \mathcal{I}}$ denote an orthonormal basis of $K$ and let $(Y,X)\in \hypo(F)(K)$. Then for an arbitrary unit vector $\gamma=\sum_{j\in\mathcal{I}}\gamma_j^*\otimes e_j\in N^*\otimes K$ we have
\begin{equation*}
\begin{split}
\gamma^*L_{B,v}(Y,X)\gamma=&\sum_{i,j\in\mathcal{I}}\gamma_j^*vv^*\gamma_ie_i^*Ye_j-\sum_{l=1}^k\gamma_j^*cB_l\gamma_ie_i^*(X_l-I)e_j\\
=&\sum_{i,j\in\mathcal{I}}\tr\left\{vv^*\gamma_ie_i^*Ye_j\gamma_j^*-\sum_{l=1}^kcB_l\gamma_ie_i^*(X_l-I)e_j\gamma_j^*\right\}\\
=&\tr\left\{vv^*\Gamma^*Y\Gamma-\sum_{l=1}^kcB_l\Gamma^*(X_l-I)\Gamma\right\}
\end{split}
\end{equation*}
where $\Gamma:N\mapsto K$ is the contraction defined as $\Gamma:=\sum_{i\in\mathcal{I}}e_i\gamma_i^*$ where convergence is in the ultraweak operator topology. Using \eqref{eq:P:separating_pencil:3} we have
\begin{equation*}
\begin{split}
\tr\left\{vv^*\Gamma^*Y\Gamma-\sum_{l=1}^kcB_l\Gamma^*(X_l-I)\Gamma\right\}&\leq \tr(\Gamma cT\Gamma^*)\\
&=\tr\left(\sum_{i,j\in\mathcal{I}}e_j\gamma_j^*cT\gamma_ie_i^*\right)\\
&=\sum_{i,j\in\mathcal{I}}\gamma_j^*cT\gamma_ie_i^*e_j\\
&=\gamma^*(cT\otimes I)\gamma.
\end{split}
\end{equation*}
Thus the linear pencil $cT-L_{B,v}$ defined by $[cT-L_{B,v}](Y,X)=cT\otimes I-vv^*\otimes Y+\sum_{i=1}^kcB_i\otimes (X_i-I)$ satisfies
\begin{equation*}
[cT-L_{B,v}](Y,X)\geq 0
\end{equation*}
for every $K\leq E$ and $(Y,X)\in \hypo(F)(K)$.


Also computing as above \eqref{eq:P:separating_pencil:4.1} becomes
\begin{equation}\label{eq:P:separating_pencil:5}
\tr\left\{\mathfrak{E}\left([cT-L_{B,v}](F(A),A)\right)\right\}=0
\end{equation}
with $\mathfrak{E}=\sum_{i,j\in\mathcal{I}}(e_i^*\otimes e_i)(e_j^*\otimes e_j)^*$.

Let $B(F,A,v)_0:=cT$ and $B(F,A,v)_i:=cB_i$. Then (b) of the assertion is satisfied. Now we turn to the proof of (a). Observe that the point $(-dI,(I,\ldots,I))$ is in $\hypo(F)$ for any scalar $d\geq 0$ since $F$ is positive, so we must have $L_{F,A,v}(-dI,(I,\ldots,I))\geq 0$ for all scalar $d\geq 0$. In other words we have $B(F,A,v)_0\otimes I\geq vv^*\otimes -dI$ for all scalar $d\geq 0$. This is only possible if $c\geq 0$ as we established earlier. Also notice that for any fixed $1\leq i\leq k$ with $Y=0$, $X_j=I$ for $j\neq i$, $X_i=dI$; the point $(Y,X)$ is in $\hypo(F)$ for any scalar $d>1$, so we must have $L_{F,A,v}(Y,X)\geq 0$ for such pairs, so it follows that $B(F,A,v)_0\otimes I+B(F,A,v)_i\otimes dI\geq 0$ always for any $d>1$ which is only possible if $B(F,A,v)_i\geq 0$. Similarly the point $(0,(dI,\ldots,dI))$ is in $\hypo(F)$ for any scalar $d>0$, hence similar consideration reveals that $B(F,A,v)_0\otimes I-(1-d)\sum_{i=1}^kB(F,A,v)_i\otimes I\geq 0$ for arbitrarily small $d>0$, thus by taking the limit $d\to 0+$ we get $\sum_{i=1}^kB(F,A,v)_i\leq B(F,A,v)_0$, finishing the proof of property (a).

It follows from \eqref{eq:P:separating_pencil:5} that $\tr\left\{\mathfrak{E}L_{B,A,v}(F(A),A)\right\}=0$ or equivalently
\begin{equation}\label{eq:P:separating_pencil:6}
\sum_{i,j\in\mathcal{I}}(e_i\otimes e_i^*)L_{B,A,v}(F(A),A)(e_j^*\otimes e_j)=0.
\end{equation}

Now assume that $c_1I\leq A_i\leq c_2I$ for all $1\leq i\leq k$ and some fixed real constants $c_2>c_1>0$ as in (c). Then by (a) and using the notation in \eqref{eq:P:separating_pencil:4.1} we have 
\begin{equation*}
\sum_{i=1}^k\tr\{B_i(A_i-I)\}\geq (c_1-1)\sum_{i=1}^k\tr(B_i)
\end{equation*}
and $1=\tr(T)\geq \sum_{i=1}^k\tr(B_i)\geq 0$. Since $A_i\leq c_2I$ for all $1\leq i\leq k$, by the operator monotonicity of $F$ we also have $v^*F(A)v\leq v^*F(c_2,\ldots,c_2)Iv=F(c_2,\ldots,c_2)$, hence using \eqref{eq:P:separating_pencil:4.1} we get
\begin{equation*}
c=\frac{v^*F(A)v}{1+\sum_{i=1}^k\tr\{B_i(A_i-I)\}}\leq \frac{F(c_2,\ldots,c_2)}{1+(c_1-1)\sum_{i=1}^k\tr(B_i)}\leq \frac{F(c_2,\ldots,c_2)}{\min(1,c_1)}
\end{equation*}
which together with $\tr(B(F,A,v)_0)=c\tr(T)$ proves (c).
\end{proof}

\section{Explicit LMI solution formula}

Proposition~\ref{P:separating_pencil} provides us a tool to find sufficiently many supporting linear pencils of hypographs of our functions so that we can reconstruct the values of the functions at each point. We will need the following result from \cite{anderson}.
\begin{theorem}[Theorem 3 cf. \cite{anderson}]\label{T:Schur_complement}
Let $A$ be a positive semi-definite linear operator on a Hilbert space and $S$ a subspace. Let the matrix of $A$ be partitioned as $A=\left[ \begin{array}{cc}
A_{11} & A_{12} \\
A_{21} & A_{22} \end{array} \right]$ with $A_{11}:S\mapsto S$, $A_{21}:S\mapsto S^{\perp}$. Then $\ran(A_{21})\subset \ran(A_{22})^{1/2}$ and there exists a bounded linear operator $C:S\mapsto S^{\perp}$ such that $A_{21}=(A_{22})^{1/2}C$ and
\begin{equation*}
A=\left[ \begin{array}{cc}
A_{11}-C^*C & 0 \\
0 & 0 \end{array} \right]+\left[ \begin{array}{cc}
C^* & 0 \\
(A_{22})^{1/2} & 0 \end{array} \right]
\left[ \begin{array}{cc}
C & (A_{22})^{1/2} \\
0 & 0 \end{array} \right].
\end{equation*}
The bounded positive semi-definite operator $\mathcal{S}(A)=A_{11}-C^*C$ is called the shorted operator or Schur complement of $A$. It satisfies $\mathcal{S}(A)\leq A$ and it is maximal among all self-adjoint operators $X:S\mapsto S$ such that $X\leq A$.
\end{theorem}

Note that the proof of the above result in \cite{anderson} is based on the following well known result.
\begin{lemma}[Douglas' lemma, cf. \cite{fillmore}]
Let $A$ and $B$ be bounded linear operators on a Hilbert space $\mathcal{H}$. Then the following statements are equivalent:
\begin{itemize}
\item[(1)]$\ran(A)\subset\ran(B)$.
\item[(2)]$AA^*\leq\lambda^2BB^*$ for some $\lambda\geq 0$.
\item[(3)]There exists a bounded linear operator $C$ such that $A=BC$.
\end{itemize}
\end{lemma}
\begin{remark}\label{R:rangeinclusion}
Note that in the proof of the Douglas' lemma the following construction is used to prove $(1)\Longrightarrow (3)$. The operator $C$ is defined as $C=B_0^{-1}A$, where $B_0$ is the restriction of $B$ to the orthogonal complement of its kernel $\mathcal{N}(B)^{\perp}$, so that $B_0^{-1}:\ran(B)\mapsto\mathcal{N}(B)^{\perp}$ is a closed linear operator, hence also $C$ is a closed linear operator from $\mathcal{H}$ into $\mathcal{N}(B)^{\perp}$. From the closed graph theorem it follows that $C$ is bounded.
\end{remark}

We will also use the following basic fact from time to time.
\begin{lemma}
Let $A\geq 0$ be a positive semi-definite operator on some Hilbert space $H$. If $v^*Av=0$ for some $v\in H$, then $Av=0$.
\end{lemma}
\begin{proof}
Since $A\geq 0$, it has a unique positive semi-definite square root $A^{1/2}$. Hence
\begin{equation*}
0=v^*Av=\|A^{1/2}v\|^2=0,
\end{equation*}
so $A^{1/2}v=0$ and it follows that $Av=0$.
\end{proof}

\begin{theorem}\label{T:reconstruct_schur}
Let $F:\mathbb{P}^k\mapsto \mathbb{P}$ be an operator monotone function. Then for each $A\in\mathbb{P}(N)^k$ with $\dim(N)<\infty$ and each unit vector $v\in N$ we have
\begin{equation}\label{eq:T:reconstruct_schur:0.1}
\begin{split}
F(A)v=&v^*{B}_{0,11}(A,v)v\otimes Iv+\sum_{i=1}^kv^*B_{i,11}(A,v)v\otimes (A_i-I)v\\
&-\left\{(v^*\otimes I)\left[{B}_{0,12}(A,v)\otimes I+\sum_{i=1}^kB_{i,12}(A,v)\otimes (A_i-I)\right]\right.\\
&\times\left[{B}_{0,22}(A,v)\otimes I+\sum_{i=1}^kB_{i,22}(A,v)\otimes (A_i-I)\right]^{-1}\\
&\left.\times\left[{B}_{0,21}(A,v)\otimes I+\sum_{i=1}^kB_{i,21}(A,v)\otimes (A_i-I)\right](v\otimes I)\right\}v
\end{split}
\end{equation}
and
\begin{equation}\label{eq:T:reconstruct_schur:0.12}
\begin{split}
&\left\{\left[{B}_{0,22}(A,v)\otimes I+\sum_{i=1}^kB_{i,22}(A,v)\otimes (A_i-I)\right]\right.\\
&\left.-\left[{B}_{0,21}(A,v)\otimes I+\sum_{i=1}^kB_{i,21}(A,v)\otimes (A_i-I)\right]\right\}(v^*\otimes v)\\
&=\sum_{j\in\mathcal{I}}\left[{B}_{0,22}(A,v)\otimes I+\sum_{i=1}^kB_{i,22}(A,v)\otimes (A_i-I)\right](e_j^*\otimes e_j),
\end{split}
\end{equation}
where $\{e_j\}_{j\in\mathcal{J}}$ is an orthonormal basis of $N$ and
\begin{equation}\label{eq:T:reconstruct_schur:0.2}
\begin{split}
{B}_{i,11}(A,v):=&vv^*{B}_{i}(A,v)vv^*,\\
{B}_{i,12}(A,v):=&vv^*{B}_{i}(A,v)(I-vv^*),\\
{B}_{i,21}(A,v):=&(I-vv^*){B}_{i}(A,v)vv^*,\\
{B}_{i,22}(A,v):=&(I-vv^*){B}_{i}(A,v)(I-vv^*)
\end{split}
\end{equation}
for all $0\leq i\leq k$ 
 and $B_{i}(A,v)=B(F,A,v)_i$.

Moreover if $c_1I\leq A_i\leq c_2I$ for all $1\leq i\leq k$ and some fixed real constants $c_2>c_1>0$, then
\begin{equation}\label{eq:T:reconstruct_schur:0.3}
\tr\{B_0(A,v)\}\leq \frac{F(c_2,\ldots,c_2)}{\min(1,c_1)}.
\end{equation}
\end{theorem}
\begin{proof}
From Proposition~\ref{P:separating_pencil} we have that
\begin{equation}\label{eq:T:reconstruct_schur:1}
vv^*\otimes F(A)\leq  B(F,A,v)_0\otimes I+\sum_{i=1}^kB(F,A,v)_i\otimes (X_i-I)
\end{equation}
and also by \eqref{eq:P:separating_pencil:6}
\begin{equation}\label{eq:T:reconstruct_schur:2}
\begin{split}
\sum_{i,j\in\mathcal{I}}(e_i\otimes e_i^*)vv^*\otimes F(A)(e_j^*\otimes e_j)=&\sum_{i,j\in\mathcal{I}}(e_i\otimes e_i^*)\left[B(F,A,v)_0\otimes I\right.\\
&\left.+\sum_{i=1}^kB(F,A,v)_i\otimes (X_i-I)\right](e_j^*\otimes e_j).
\end{split}
\end{equation}
By (a) of Proposition~\ref{P:separating_pencil} we can apply the Schur complement Theorem~\ref{T:Schur_complement} to \eqref{eq:T:reconstruct_schur:1}, pivoting on the subspace $v^*\otimes N$ of $N^*\otimes N$ to get
\begin{equation}\label{eq:T:reconstruct_schur:3}
\begin{split}
vv^*\otimes F(A)\leq &{B}_{0,11}(A,v)\otimes I+\sum_{i=1}^kB_{i,11}(A,v)\otimes (A_i-I)\\
&-\left[{B}_{0,12}(A,v)\otimes I+\sum_{i=1}^kB_{i,12}(A,v)\otimes (A_i-I)\right]\\
&\times\left[{B}_{0,22}(A,v)\otimes I+\sum_{i=1}^kB_{i,22}(A,v)\otimes (A_i-I)\right]^{-1}\\
&\times\left[{B}_{0,21}(A,v)\otimes I+\sum_{i=1}^kB_{i,21}(A,v)\otimes (A_i-I)\right]
\end{split}
\end{equation}
where the coefficients 
 ${B}_{i,xy}(A,v)$ are defined by \eqref{eq:T:reconstruct_schur:0.2} and $B_{i}(A,v):=B(F,A,v)_i$. The inversion in \ref{eq:T:reconstruct_schur:3} is justified by (a) in Proposition~\ref{P:separating_pencil} alone or by the range inclusion result in Theorem~\ref{T:Schur_complement} and that the operators have finite rank. Notice that
\begin{equation*}
\begin{split}
(vv^*\otimes I)\left[\sum_{i,j\in\mathcal{I}}(e_j^*\otimes e_j)(e_i\otimes e_i^*)\right](vv^*\otimes I)=&\sum_{i,j\in\mathcal{I}}(e_j^*vv^*\otimes e_j)(vv^*e_i\otimes e_i^*)\\
=&\sum_{i,j\in\mathcal{I}}(v^*\otimes e_je_j^*v)(v\otimes v^*e_ie_i^*),\\
=&(v^*\otimes v)(v\otimes v^*)
\end{split}
\end{equation*}
hence from \eqref{eq:T:reconstruct_schur:1}, \eqref{eq:T:reconstruct_schur:2} and \eqref{eq:T:reconstruct_schur:3} we get
\begin{equation}\label{eq:T:reconstruct_schur:4}
\begin{split}
v^*F(A)v=&v^*{B}_{0,11}(A,v)v\otimes v^*Iv+\sum_{i=1}^kv^*B_{i,11}(A,v)v\otimes v^*(A_i-I)v\\
&-(v\otimes v^*)\left[{B}_{0,12}(A,v)\otimes I+\sum_{i=1}^kB_{i,12}(A,v)\otimes (A_i-I)\right]\\
&\times\left[{B}_{0,22}(A,v)\otimes I+\sum_{i=1}^kB_{i,22}(A,v)\otimes (A_i-I)\right]^{-1}\\
&\times\left[{B}_{0,21}(A,v)\otimes I+\sum_{i=1}^kB_{i,21}(A,v)\otimes (A_i-I)\right](v^*\otimes v)+e
\end{split}
\end{equation}
where
\begin{equation*}
\begin{split}
e:=&\left\{(v\otimes v^*)\left[{B}_{0,12}(A,v)\otimes I+\sum_{i=1}^kB_{i,12}(A,v)\otimes (A_i-I)\right]\right.\\
&\times\left[{B}_{0,22}(A,v)\otimes I+\sum_{i=1}^kB_{i,22}(A,v)\otimes (A_i-I)\right]^{-1/2}\\
&\left.+\left[\sum_{j\in\mathcal{I}}(e_j\otimes e_j^*)-(v\otimes v^*)\right]\left[{B}_{0,22}(A,v)\otimes I+\sum_{i=1}^kB_{i,22}(A,v)\otimes (A_i-I)\right]^{1/2}\right\}\\
&\left\{\left[{B}_{0,22}(A,v)\otimes I+\sum_{i=1}^kB_{i,22}(A,v)\otimes (A_i-I)\right]^{-1/2}\right.\\
&\times\left[{B}_{0,21}(A,v)\otimes I+\sum_{i=1}^kB_{i,21}(A,v)\otimes (A_i-I)\right](v^*\otimes v)\\
&\left.+\left[{B}_{0,22}(A,v)\otimes I+\sum_{i=1}^kB_{i,22}(A,v)\otimes (A_i-I)\right]^{1/2}\left[\sum_{j\in\mathcal{I}}(e_j^*\otimes e_j)-(v^*\otimes v)\right]\right\}.
\end{split}
\end{equation*}
Notice that $e$ is of the form $e=x^*x$, where
\begin{equation}\label{eq:T:reconstruct_schur:5}
\begin{split}
x=&D^{-1/2}\left[{B}_{0,21}(A,v)\otimes I+\sum_{i=1}^kB_{i,21}(A,v)\otimes (A_i-I)\right](v^*\otimes v)\\
&+D^{1/2}\left[\sum_{j\in\mathcal{I}}(e_j^*\otimes e_j)-(v^*\otimes v)\right],\\
D=&\left[{B}_{0,22}(A,v)\otimes I+\sum_{i=1}^kB_{i,22}(A,v)\otimes (A_i-I)\right].
\end{split}
\end{equation}
It follows that $e\geq 0$. The linear map $l:\mathcal{B}(N^*)\otimes \mathcal{B}(N)\mapsto \mathbb{C}\otimes \mathcal{B}(N)\cong \mathcal{B}(N)$ defined on simple tensors as $l(X\otimes Y):=v^*Xv\otimes Y$ is (completely) positive, hence order preserving. So applying $l$ to \eqref{eq:T:reconstruct_schur:3} we get
\begin{equation}\label{eq:T:reconstruct_schur:6}
\begin{split}
F(A)\leq&v^*{B}_{0,11}(A,v)v\otimes I+\sum_{i=1}^kv^*B_{i,11}(A,v)v\otimes (A_i-I)\\
&-(v^*\otimes I)\left[{B}_{0,12}(A,v)\otimes I+\sum_{i=1}^kB_{i,12}(A,v)\otimes (A_i-I)\right]\\
&\times\left[{B}_{0,22}(A,v)\otimes I+\sum_{i=1}^kB_{i,22}(A,v)\otimes (A_i-I)\right]^{-1}\\
&\times\left[{B}_{0,21}(A,v)\otimes I+\sum_{i=1}^kB_{i,21}(A,v)\otimes (A_i-I)\right](v\otimes I).
\end{split}
\end{equation}
Now since we have \eqref{eq:T:reconstruct_schur:6} and also \eqref{eq:T:reconstruct_schur:4} with $e\geq 0$, we must have $e=0$ in \eqref{eq:T:reconstruct_schur:4} which yields $x=0$ in \eqref{eq:T:reconstruct_schur:5} hence \eqref{eq:T:reconstruct_schur:0.12} and also \eqref{eq:T:reconstruct_schur:0.1}.

Property (c) in Proposition~\ref{P:separating_pencil} yields $\tr\{B_0(A,v)\}\leq \frac{F(c_2,\ldots,c_2)}{\min(1,c_1)}$.
\end{proof}

\begin{definition}[Natural map]
A graded map ${F}:\mathbb{S}(K)^k\times K\mapsto K$ for each Hilbert space $K$ is called a \emph{natural map} if it preserves direct sums, i.e.
\begin{equation*}
{F}(X\oplus Y,v\oplus w)={F}(X,v)\oplus {F}(Y,w)
\end{equation*}
for $X\in\mathbb{S}(K_1)^k$, $v\in K_1$ and $Y\in\mathbb{S}(K_2)^k$, $w\in K_2$, for Hilbert spaces $K_1,K_2$.
\end{definition}

For an NC function $F:\mathbb{S}^k\mapsto \mathbb{S}$ we define the natural map $\overline{F}:\mathbb{S}(K)^k\times K\mapsto K$ for any Hilbert space $K$ by
\begin{equation*}
\overline{F}(X,v):=F(X)v
\end{equation*}
for $X\in\mathbb{S}(K)^k$ and $v\in K$. Indeed $\overline{F}$ is natural since
\begin{equation*}
\overline{F}(X\oplus Y,v\oplus w)=F(X)v\oplus F(Y)w=\overline{F}(X,v)\oplus \overline{F}(Y,w)
\end{equation*}
for $X\in\mathbb{S}(K_1)^k$, $v\in K_1$ and $Y\in\mathbb{S}(K_2)^k$, $w\in K_2$ for Hilbert spaces $K_1,K_2$.
Also notice that the function
\begin{equation*}
\begin{split}
F(X):=&v^*{B}_{0,11}v\otimes I+\sum_{i=1}^kv^*B_{i,11}v\otimes (X_i-I)\\
&-(v^*\otimes I)\left[{B}_{0,12}\otimes I+\sum_{i=1}^kB_{i,12}\otimes (X_i-I)\right]\\
&\times\left[{B}_{0,22}\otimes I+\sum_{i=1}^kB_{i,22}\otimes (X_i-I)\right]^{-1}\\
&\times\left[{B}_{0,21}\otimes I+\sum_{i=1}^kB_{i,21}\otimes (X_i-I)\right](v\otimes I)
\end{split}
\end{equation*}
for fixed $v\in E$ and linear operators $B_{i,xy}$ such that the products above are well defined is an NC function, for example like in \eqref{eq:T:reconstruct_schur:0.1} of Theorem~\ref{T:reconstruct_schur}. Hence it defines a natural map in the same way by
\begin{equation*}
\overline{F}(X,w):=F(X)w
\end{equation*}
for $X\in\mathbb{P}(K)^k$ and $w\in K$.

Let $S(E):=\{v\in E:\|v\|=1\}$ denote the unit sphere of the Hilbert space $E$. In what follows, we will construct an auxiliary Hilbert space on which we will apply Theorem~\ref{T:reconstruct_schur}.

For fixed real constants $c_2>c_1>0$, let
\begin{equation*}
\begin{split}
\mathbb{P}_{c_1,c_2}(E)&:=\{X\in\mathbb{P}(E):c_1I\leq X\leq c_2I\},\\
\Omega_{c_1,c_2}&:=\mathbb{P}_{c_1,c_2}(E)^k\times S(E)
\end{split}
\end{equation*}
and let
\begin{equation*}
\mathcal{H}:=\bigoplus_{\dim(E)<\infty}\bigoplus_{\omega\in\Omega_{c_1,c_2}}E.
\end{equation*}
We equip $\mathcal{H}$ with the inner product
\begin{equation*}
x^*y:=\sum_{\dim(E)<\infty}\sum_{\omega\in\Omega_{c_1,c_2}}x(\omega)^*y(\omega)
\end{equation*}
for $x,y\in \mathcal{H}$ and we denote again by $\mathcal{H}$ the Hilbert space completion with respect to this inner product.


\begin{definition}\label{D:Psi}
Let $F:\mathbb{P}^k\mapsto \mathbb{P}$ be an operator monotone function. Now let
\begin{equation*}
\begin{split}
\Psi_{F}(X):=&\bigoplus_{\dim(E)<\infty}\bigoplus_{(A,v)\in\Omega_{c_1,c_2}}\left\{{B}_{0,11}(A,v)\otimes I+\sum_{i=1}^kB_{i,11}(A,v)\otimes (X_i-I)\right.\\
&-\left[{B}_{0,12}(A,v)\otimes I+\sum_{i=1}^kB_{i,12}(A,v)\otimes (X_i-I)\right]\\
&\times\left[{B}_{0,22}(A,v)\otimes I+\sum_{i=1}^kB_{i,22}(A,v)\otimes (X_i-I)\right]^{-1}\\
&\left.\times\left[{B}_{0,21}(A,v)\otimes I+\sum_{i=1}^kB_{i,21}(A,v)\otimes (X_i-I)\right]\right\},
\end{split}
\end{equation*}
where the coefficients 
 ${B}_{i,xy}(A,v)$ for $0\leq i\leq k$ and $x,y\in\{1,2\}$ are as in \eqref{eq:T:reconstruct_schur:0.1} of Theorem~\ref{T:reconstruct_schur}. Then $\Psi_{F}(X)$ is a linear operator on $\mathcal{H}\otimes E$ which is also bounded by \eqref{eq:T:reconstruct_schur:0.3} and Theorem~\ref{T:Schur_complement} for $X\in\mathbb{P}(E)^k$. We also have
\begin{equation*}
\begin{split}
\Psi_{F}(X)=&{B}_{0,11}\otimes I+\sum_{i=1}^kB_{i,11}\otimes (X_i-I)\\
&-\left[{B}_{0,12}\otimes I+\sum_{i=1}^kB_{i,12}\otimes (X_i-I)\right]\\
&\times\left[{B}_{0,22}\otimes I+\sum_{i=1}^kB_{i,22}\otimes (X_i-I)\right]^{-1}\\
&\times\left[{B}_{0,21}\otimes I+\sum_{i=1}^kB_{i,21}\otimes (X_i-I)\right]
\end{split}
\end{equation*}
where
\begin{equation*}
{B}_{i,xy}:=\bigoplus_{\dim(E)<\infty}\bigoplus_{(A,v)\in\Omega_{c_1,c_2}}{B}_{i,xy}(A,v)
\end{equation*}
for $0\leq i\leq k$ and $x,y\in\{1,2\}$.
\end{definition}

\begin{lemma}\label{L:reconstruct_schur_multiple}
Let $F:\mathbb{P}^k\mapsto \mathbb{P}$ be an operator monotone function and let $\dim(E)<\infty$. Let $A_j\in\mathbb{P}_{c_1,c_2}(E)^k$ and $v_j\in S(E)$ for $j\in \mathcal{J}$ for some finite index set $\mathcal{J}$. Then there exists a $w\in S(\mathcal{H})$ such that
\begin{equation}\label{eq:L:reconstruct_schur_multiple:0.1}
F(A_j)v_j=(w^*\otimes I)\Psi_{F}(A_j)(w\otimes I)v_j
\end{equation}
for all $j\in \mathcal{J}$.
\end{lemma}
\begin{proof}
Let $A:=\bigoplus_{j\in \mathcal{J}}A_j$ and $v:=\bigoplus_{j\in \mathcal{J}}\frac{1}{\sqrt{|\mathcal{J}|}}v_j$. Then by Theorem~\ref{T:reconstruct_schur} and the definition of $\Psi_F$, we have
\begin{equation*}
\begin{split}
F(A)v=&v^*{B}_{0,11}(A,v)v\otimes Iv+\sum_{i=1}^kv^*B_{i,11}(A,v)v\otimes (A_i-I)v\\
&-\left\{(v^*\otimes I)\left[{B}_{0,12}(A,v)\otimes I+\sum_{i=1}^kB_{i,12}(A,v)\otimes (A_i-I)\right]\right.\\
&\times\left[{B}_{0,22}(A,v)\otimes I+\sum_{i=1}^kB_{i,22}(A,v)\otimes (A_i-I)\right]^{-1}\\
&\left.\times\left[{B}_{0,21}(A,v)\otimes I+\sum_{i=1}^kB_{i,21}(A,v)\otimes (A_i-I)\right](v\otimes I)\right\}v\\
&=(w^*\otimes I)\Psi_{F}(A)(w\otimes I)v,
\end{split}
\end{equation*}
where $w:=(\cdots\oplus 0 \oplus v \oplus 0\oplus\cdots)\in S(\mathcal{H})$, and the nonzero $v$ is at the appropriate coordinate such that the above holds, according to the definition of $\Psi_{F}$. Since both $(A,v)\mapsto F(A)v$ and $(A,v)\mapsto (w^*\otimes I)\Psi_{F}(A)(w\otimes I)v$ are natural maps, we have \eqref{eq:L:reconstruct_schur_multiple:0.1}.
\end{proof}

Let $\mathcal{B}^+_1(\mathcal{H})^{*}$ denote the state space of $\mathcal{B}(\mathcal{H})$ and $\mathcal{B}^+_1(\mathcal{H})_{*}$ its normal part. Note that positive linear functionals on unital $C^*$-algebras attain their norm at the unit, hence $\mathcal{B}_1^+(\mathcal{H})^{*}$ is convex, weak-$*$ compact by Banach-Alaoglu.

\begin{theorem}\label{T:Loewner_formula}
Let $F:\mathbb{P}^k\mapsto \mathbb{P}$ be an operator monotone function. Then there exists a $w\in \mathcal{B}^+_1(\mathcal{H})^{*}$ such that for all $\dim(E)<\infty$ and $X\in\mathbb{P}_{c_1,c_2}(E)^k$ we have
\begin{equation}\label{eq:T:Loewner_formula:0.1}
\begin{split}
F(X)=&(w\otimes I)(\Psi_{F}(X))\\
=&w({B}_{0,11})\otimes I+\sum_{i=1}^kw(B_{i,11})\otimes (X_i-I)\\
&-(w\otimes I)\left\{\left[{B}_{0,12}\otimes I+\sum_{i=1}^kB_{i,12}\otimes (X_i-I)\right]\right.\\
&\times\left[{B}_{0,22}\otimes I+\sum_{i=1}^kB_{i,22}\otimes (X_i-I)\right]^{-1}\\
&\left.\times\left[{B}_{0,21}\otimes I+\sum_{i=1}^kB_{i,21}\otimes (X_i-I)\right]\right\}.
\end{split}
\end{equation}
\end{theorem}
\begin{proof}
For $A\in\mathbb{P}_{c_1,c_2}(E)^k$ and $v\in E$ by Lemma~\ref{L:reconstruct_schur_multiple} the set
\begin{equation*}
L_{(A,v)}:=\{w\in \mathcal{B}^+_1(\mathcal{H})^{*}:(w\otimes I)(\Psi_{F}(A))v=F(A)v\}
\end{equation*}
is nonempty, moreover it is easy to check that it is a closed subset of $\mathcal{B}^+_1(\mathcal{H})^{*}$ in the weak-$*$ topology of $\mathcal{B}(\mathcal{H})^{*}$ hence compact. Indeed, the latter is the consequence of the identification $\mathcal{B}(\mathcal{H}\otimes E)=\mathcal{B}(\mathcal{H}^{\dim(E)})\simeq M_{\dim(E)}(\mathcal{B}(\mathcal{H}))$, where $M_{\dim(E)}(\mathcal{B}(\mathcal{H}))$ denotes the $\dim(E)$-by-$\dim(E)$ matrices with $\mathcal{B}(\mathcal{H})$ valued entries.

Let $\bold{L}$ denote the collection $\{L_{(A,v)}:A\in\mathbb{P}_{c_1,c_2}(E)^k,v\in E\}$ of subsets of $\mathcal{B}^+_1(\mathcal{H})^{*}$. Any finite sub-collection from $\bold{L}$ has the form $\{L_{(A,v)}:(A,v)\in\Omega_{c_1,c_2}^n\}$ for some integer $n$, so by Lemma~\ref{L:reconstruct_schur_multiple} has nonempty intersection. This means that $\bold{L}$ has the finite intersection property. So the weak-$*$ compactness of $\mathcal{B}^+_1(\mathcal{H})^{*}$ implies that there is a $w\in \mathcal{B}^+_1(\mathcal{H})^{*}$ which is in every $L_{(A,v)}$, proving the assertion.
\end{proof}

\section{Analytic properties of resolvents}

\begin{definition}
The \emph{imaginary part} of a bounded linear operator $A$ is defined as
\begin{equation*}
\Im A:=\frac{A-A^*}{2i},
\end{equation*}
and its \emph{real part} is
\begin{equation*}
\Re A:=\frac{A+A^*}{2}.
\end{equation*}
\end{definition}

\begin{proposition}\label{P:Schur_imgainary_invariant}
Let
\begin{equation*}
A:=\left[ \begin{array}{cc}
A_{11} & A_{12} \\
A_{21} & A_{22} \end{array} \right]
\end{equation*}
be a block operator matrix with $\Im A\geq 0$. Then if its Schur complement
\begin{equation*}
\mathcal{S}(A)=A_{22}-A_{21}A_{11}^{-1}A_{12}
\end{equation*}
exists, it also satisfies $\Im \mathcal{S}(A)\geq 0$.

Similarly if $\Re A\geq 0$ and if its Schur complement exists, it also satisfies $\Re \mathcal{S}(A)\geq 0$.
\end{proposition}
\begin{proof}
By assumption we have
\begin{equation*}
0\leq \Im \left\{(x_1\oplus x_2)^*\left[ \begin{array}{cc}
A_{11} & A_{12} \\
A_{21} & A_{22} \end{array} \right](x_1\oplus x_2)\right\}
\end{equation*}
for any vector $(x_1\oplus x_2)$. Hence choosing $x_1=-A_{11}^{-1}A_{12}x_2$ we get that
\begin{equation*}
\begin{split}
0&\leq \Im \left\{(x_1\oplus x_2)^*\left[ \begin{array}{cc}
A_{11} & A_{12} \\
A_{21} & A_{22} \end{array} \right](x_1\oplus x_2)\right\}\\
&=\Im \left\{(x_1\oplus x_2)^*\left[ \begin{array}{cc}
0\\
A_{22}-A_{21}A_{11}^{-1}A_{12}\end{array} \right]\right\}\\
&=\Im \{x_2^*\mathcal{S}(A)x_2\}
\end{split}
\end{equation*}
proving the first part of the assertion. The second part covering the real parts is proved similarly.
\end{proof}

The following result for the Schur complement is also well known and its proof can be found for example in \cite{anderson}.
\begin{proposition}\label{P:Schur_concave_monotone}
Let
\begin{equation*}
A:=\left[ \begin{array}{cc}
A_{11} & A_{12} \\
A_{21} & A_{22} \end{array} \right]\text{ and }
B:=\left[ \begin{array}{cc}
B_{11} & B_{12} \\
B_{21} & B_{22} \end{array} \right]
\end{equation*}
be conformally partitioned, positive semi-definite block operator matrices. Then if $A\leq B$ then also $\mathcal{S}(A)\leq \mathcal{S}(B)$, i.e. the Schur complement is operator monotone.

Moreover $\mathcal{S}(\cdot)$ is also operator concave, i.e. $(1-\lambda)\mathcal{S}(A)+\lambda \mathcal{S}(B)\leq \mathcal{S}((1-\lambda)A+\lambda B)$ for all $\lambda\in [0,1]$.
\end{proposition}

\begin{definition}[Sectorial operator]
For $A\in\mathcal{B}(E)$ let $W(A):=\{x^*Ax:x\in E,\|x\|=1\}$ denote the numerical range of $A$. We say that $A$ is sectorial if $W(A)\subseteq S_\alpha$, where $S_\alpha:=\{z\in\mathbb{C}:\Re(z)>0, |\Im(z)|\leq Re(z)\tan \alpha\}$ for some $\alpha\in[0,\pi/2)$.
\end{definition}

For a matrix $X\in\mathcal{B}(N)$ where $\dim(N)<\infty$ let $\{\sigma_j(X)\}_{1\leq j\leq \dim(N)}$ denote the ordered decreasing sequence of its singular values.

\begin{proposition}[Theorem 1.1 \cite{drury}]\label{P:Drury}
Let $A\in\mathcal{B}(N)$ with $\dim(N)<\infty$ and $W(A)\subseteq S_{\alpha}$ be partitioned as
\begin{equation*}
A:=\left[ \begin{array}{cc}
A_{11} & A_{12} \\
A_{21} & A_{22} \end{array} \right].
\end{equation*}
Let $\mathcal{S}(A)=A_{22}-A_{21}A_{11}^{-1}A_{12}$. Then
\begin{equation*}
\sigma_j(\mathcal{S}(A))\leq \sec^2(\alpha)\sigma_j(A_{22}).
\end{equation*}
In particular under the assumptions above including also $\dim(N)=\infty$ and that $A$ has a bounded inverse, we have
\begin{equation}\label{eq:P:Drury}
\|\mathcal{S}(A)\|\leq \sec^2(\alpha)\|A\|.
\end{equation}
\end{proposition}
\begin{proof}
The proof of this under the assumption $\dim(N)<\infty$ can be found in \cite{drury}. Invertibility follows from $W(A)\subseteq S_{\alpha}$ and the result in \cite{drury} covers the inequality for the singular values and from that \eqref{eq:P:Drury} follows for the finite dimensional case, since $\|A_{22}\|\leq \|A\|$.

If $\dim(N)=\infty$, we assume that the block operator matrix $A$ is partitioned according to the orthogonal decomposition $N=S\oplus S^{\perp}$, where $S$ is a closed subspace of $N$. We can approximate $A$ in the strong operator topology by a net of finite rank operators $P_{\gamma}AP_{\gamma}$, where $P_{\gamma}$ is the directed set of projections with ranges running over all finite dimensional subspaces of $N$ partially ordered under subspace inclusion. Indeed, $A_\gamma:=P_{\gamma}AP_{\gamma}$ converges to $A$ in the strong operator topology, since
\begin{equation*}
\begin{split}
\|P_{\gamma}AP_{\gamma}x-Ax\|&\leq \|P_{\gamma}AP_{\gamma}x-P_{\gamma}Ax\|+\|P_{\gamma}Ax-Ax\|\\
&\leq \|P_{\gamma}A\|\|P_{\gamma}x-x\|+\|P_{\gamma}Ax-Ax\|\\
&\leq \|A\|\|P_{\gamma}x-x\|+\|P_{\gamma}Ax-Ax\|
\end{split}
\end{equation*}
and $\|P_{\gamma}y-y\|\to 0$ for all $y\in N$ by construction. Notice that $A_\gamma=P_{\gamma}AP_{\gamma}$ is also sectorial on the closed subspace $P_{\gamma}N\leq N$ and existence of the Schur complement $\mathcal{S}(A_\gamma)$ is justified by Theorem 1 in \cite{ando0}. 

We claim that $(A_\gamma)_{11}^{-1}\to A_{11}^{-1}$ in the strong operator topology. To see this let $P_{S}$ denote the orthogonal projection such that $A_{11}=P_{S}AP_{S}$. Then $(A_\gamma)_{11}=P_{S}P_{\gamma}AP_{\gamma}P_{S}\subseteq S_{\alpha}$ on $P_{\gamma}P_{S}N\leq N$. Then for any $z\in P_{\gamma}P_{S}N$ with $\|z\|=1$ we have
\begin{equation}\label{eq:P:Drury:3}
|z^*(A_\gamma)_{11}z|\leq \|z\|\|(A_\gamma)_{11}z\|=\|(A_\gamma)_{11}z\|.
\end{equation}
Since $W((A_\gamma)_{11})$ is a closed and compact subset of the open sector $S_{\alpha}$, we have that $|W((A_\gamma)_{11})|$ has a greatest lower bound $c>0$ that is attained in $|W((A_\gamma)_{11})|$. Hence from \eqref{eq:P:Drury:3} we have that $(A_\gamma)_{11}$ is lower bounded, i.e. $\|(A_\gamma)_{11}z\|\geq c\|z\|$ for any $z\in P_{\gamma}P_{S}N$, hence $(A_\gamma)_{11}$ is invertible on $P_{\gamma}P_{S}N\leq N$ and we have $(A_\gamma)_{11}^{-1}=((A_\gamma)_{11}^*(A_\gamma)_{11})^{-1}(A_\gamma)_{11}^*$. Similarly since $A_{11}$ is invertible on $P_{S}N\leq N$ it is lower bounded and $A_{11}^{-1}=(A_{11}^*A_{11})^{-1}A_{11}^*$. We also have that
\begin{equation*}
\begin{split}
\|(A_\gamma)_{11}^*(A_\gamma)_{11}\|&=\|P_{S}P_{\gamma}A^*P_{\gamma}P_{S}P_{S}P_{\gamma}AP_{\gamma}P_{S}\|\\
&\leq \|A^*P_{\gamma}P_{S}P_{S}P_{\gamma}A\|\\
&\leq \|A^*A\|
\end{split}
\end{equation*}
where to obtain the last inequality we used that $I\geq P_{\gamma}P_{S}P_{S}P_{\gamma}$. Similarly
\begin{equation*}
\|A_{11}^*A_{11}\|\leq \|A^*A\|,
\end{equation*}
hence it follows that
\begin{equation*}
\begin{split}
\frac{1}{2}I_{P_SN}&\geq \frac{1}{2\|A^*A\|}A_{11}^*A_{11}=:X\\
\frac{1}{2}I_{P_\gamma P_SN}&\geq \frac{1}{2\|A^*A\|}(A_\gamma)_{11}^*(A_\gamma)_{11}=:X_\gamma,
\end{split}
\end{equation*}
moreover we have that $X$ and $X_\gamma$ are lower bounded. Hence we have for some $1>t,t_\gamma>0$ that
\begin{equation*}
\begin{split}
\|I_{P_SN}-X\|&<t\\
\|I_{P_\gamma P_SN}-X_\gamma\|&<t_{\gamma},
\end{split}
\end{equation*}
so it follows that
\begin{equation}\label{eq:P:Drury:1}
\begin{split}
X^{-1}&=\left[I_{P_SN}-\left(I_{P_SN}-X\right)\right]^{-1}\\
&=\sum_{l=0}^\infty\left(I_{P_SN}-X\right)^l,\\
X_{\gamma}^{-1}&=\left[I_{P_\gamma P_SN}-\left(I_{P_\gamma P_SN}-X_{\gamma}\right)\right]^{-1}\\
&=\sum_{l=0}^\infty\left(I_{P_\gamma P_SN}-X_{\gamma}\right)^l\end{split}
\end{equation}
where $I_{P_SN}$ and $I_{P_\gamma P_SN}$ are the identity operators on the respective Hilbert subspaces $P_SN$ and $P_\gamma P_SN$. Since operator multiplication is jointly strong operator continuous on bounded sets, we have $X_{\gamma}\to X$ and $I_{P_\gamma P_SN}\to I_{P_SN}$, hence 
\begin{equation}\label{eq:P:Drury:2}
\left(I_{P_\gamma P_SN}-X_{\gamma}\right)^l\to \left(I_{P_SN}-X\right)^l
\end{equation}
in the strong operator topology for all $l\geq 0$, since also $P_\gamma\to I$ in the strong operator topology. The sums in \eqref{eq:P:Drury:1} are uniformly convergent in the norm topology, hence also in the strong operator topology, so it follows by \eqref{eq:P:Drury:2} that $X_\gamma^{-1}\to X^{-1}$ and then $((A_\gamma)_{11}^*(A_\gamma)_{11})^{-1}\to (A_{11}^*A_{11})^{-1}$ in the strong operator topology. We have that $(A_\gamma)_{11}^*=P_SP_{\gamma}A^*P_{\gamma}P_S\to P_SA^*P_S=A_{11}^*$ in the strong operator topology, then by the joint strong operator continuity of operator multiplication on bounded sets we get that $(A_\gamma)_{11}^{-1}\to A_{11}^{-1}$ in the strong operator topology as claimed.

Using the claim that $(A_\gamma)_{11}^{-1}\to A_{11}^{-1}$ in the strong operator topology and the joint strong operator continuity of operator multiplication on bounded sets we obtain
\begin{equation*}
\mathcal{S}(A_\gamma)\to \mathcal{S}(A)
\end{equation*}
where $\mathcal{S}(A)=A_{22}-A_{21}A_{11}^{-1}A_{12}$ and $\mathcal{S}(A_\gamma)=(A_\gamma)_{22}-(A_\gamma)_{21}(A_\gamma)_{11}^{-1}(A_\gamma)_{12}$.
By \eqref{eq:P:Drury} $\|\mathcal{S}(A_\gamma)\|\leq \sec^2(\alpha)\|A_\gamma\|$ and since $\|A_\gamma\|\leq \|A\|$ we have
\begin{equation*}
\|\mathcal{S}(A_\gamma)\|\leq \sec^2(\alpha)\|A\|.
\end{equation*}
Since $\mathcal{S}(A_\gamma)z\to \mathcal{S}(A)z$ for any vector $z$ we have $\|\mathcal{S}(A_\gamma)z\|\to \|\mathcal{S}(A)z\|$, hence
\begin{equation*}
\|\mathcal{S}(A)z\|\leq \sec^2(\alpha)\|A\|
\end{equation*}
for any $\|z\|=1$ and we obtain \eqref{eq:P:Drury} in the infinite dimensional case as well.
\end{proof}

\begin{lemma}\label{L:sectorial}
Let $L_{B}(X)=\sum_{i=1}^kB_i\otimes X_i$ be a linear matrix pencil over the Hilbert space $\mathcal{K}\otimes E$ such that $B_i\geq 0$. If all $X_i$ are sectorial, then $L_{B}(X)$ is also sectorial on $\mathcal{N}^{\perp}(L_{B}(1))\otimes E$, where $\mathcal{N}^{\perp}(L_{B}(1))$ denotes the closure of the complement of the kernel of $L_{B}(1)=\sum_{i=1}^kB_i\otimes 1=\sum_{i=1}^kB_i$.

Suppose additionally to the above that $\sum_{i=1}^k B_i$ is lower bounded on $\mathcal{N}^{\perp}(L_{B}(1))$. Then if all $\Re X_i$ is lower bounded, then $\Re L_{B}(X)$ is also lower bounded on $\mathcal{N}^{\perp}(L_{B}(1))\otimes E$.
\end{lemma}
\begin{proof}
We begin with a simple observation. If $A\in\mathcal{B}(N)$ is a sectorial operator on a Hilbert space $N$, then for any $X\in\mathcal{B}(N)$ the operator $X^*AX$ is also sectorial.

Now for $\gamma=\sum_{j\in\mathcal{I}}e_j^*\otimes\gamma_j\in\mathcal{N}^{\perp}(L_{B}(1))\otimes E$ we have
\begin{equation*}
\begin{split}
\gamma^*L_{B}(X)\gamma&=\sum_{i=1}^k\tr\{B_i\Gamma^*X_i\Gamma\}\\
&=\sum_{i=1}^k\tr\{B_i^{1/2}\Gamma^*X_i\Gamma B_i^{1/2}\}\\
&=\sum_{i=1}^k\sum_{j\in\mathcal{I}}{e}_{j}^*B_i^{1/2}\Gamma^*X_i\Gamma B_i^{1/2}{e}_j
\end{split}
\end{equation*}
where $\Gamma=\sum_{j\in\mathcal{I}}\gamma_je_j^*$ and convergence is in the ultraweak operator topology, $\tr\Gamma^*\Gamma=\|\gamma\|$ and $\{{e}_j\}_{j\in\mathcal{I}}$ denotes an orthonormal basis of $\mathcal{N}^{\perp}(L_{B}(1))$. From the above it follows that $W(L_{B}(X))\subseteq S_{\alpha}$ if $X_i\subseteq S_{\alpha}$ for all $1\leq i\leq k$ and a fixed $\alpha\in[0,\pi/2)$.

Now we prove the second part of the assertion. Since each $X_i$ is sectorial, we have that $\Re X_i\geq 0$. Moreover each $\Re X_i$ is lower bounded, hence there exists a real number $\epsilon >0$ such that for all $1\leq i\leq k$ we have $\Re X_i \geq \epsilon I$. Then we have that
\begin{equation}\label{eq:L:sectorial:1}
\Re L_{B}(X)=\sum_{i=1}^k B_i\otimes \Re X_i\geq \sum_{i=1}^k B_i\otimes \epsilon I.
\end{equation}
Since $\sum_{i=1}^k B_i$ is lower bounded on $\mathcal{N}^{\perp}(L_{B}(1))$, it follows from \eqref{eq:L:sectorial:1} that $\Re L_{B}(X)$ is lower bounded on $\mathcal{N}^{\perp}(L_{B}(1))\otimes E$ as well.
\end{proof}

\begin{definition}[Operator poly-halfspaces]
For a Hilbert space $E$ the \emph{upper operator poly-halfspace} is defined as
\begin{equation*}
\Pi^k:=\{X\in\mathcal{B}(E)^k:\Im X_i>0, 1\leq i\leq k\},
\end{equation*}
while the \emph{right operator poly-halfspace} as
\begin{equation*}
\Sigma^k:=\{X\in\mathcal{B}(E)^k:\Re X_i>0, 1\leq i\leq k\}.
\end{equation*}
We also use the notation $\Pi:=\Pi^1$ and $\Sigma:=\Sigma^1$.
\end{definition}

\begin{proposition}\label{P:DruryPsi}
Let $L_{B}(X)=B_0\otimes I+\sum_{i=1}^kB_i\otimes X_i$ be a linear matrix pencil over the Hilbert space $\mathcal{K}\otimes E$ such that $B_i\geq 0$. Assume that both $B_0$ and $\sum_{i=1}^k B_i$ are lower bounded on $\mathcal{N}^{\perp}(B_0)$ and $\mathcal{N}^{\perp}(\sum_{i=1}^k B_i)$ respectively. Let $S$ be a Hilbert subspace of $\mathcal{K}$ and let $P_S\in\mathcal{B}(\mathcal{K})$ denote the orthogonal projection onto $S$. Then for each $X\in\Sigma^k$ or $X\in\Pi^k$ the Schur complement
\begin{equation*}
\mathcal{S}(L_{B}(X))=PL_{B}(X)P-PL_{B}(X)P^{\perp}\left[P^{\perp}L_{B}(X)P^{\perp}\right]^{-1}P^{\perp}L_{B}(X)P
\end{equation*}
exists, where $P=P_S\otimes I$, $P^{\perp}=(I-P_S)\otimes I$, moreover there exists an $\alpha\in[0,\pi/2)$ depending on $X$ but independent of $B_i$ such that
\begin{equation}\label{eq:P:DruryPsi:1}
\|\mathcal{S}(L_{B}(X))\|\leq \sec^2\alpha\|L_{B}(X)\|.
\end{equation}
\end{proposition}
\begin{proof}
First assume that $X\in\Sigma^k$. Then there exists an $\epsilon>0$ such that $\Re X_i\geq \epsilon I$. Then
\begin{equation*}
\begin{split}
\Re L_{B}(X)&=B_0\otimes I+\sum_{i=1}^kB_i\otimes \Re X_i\\
&\geq B_0\otimes I+\sum_{i=1}^kB_i\otimes \epsilon I=\left(B_0+\epsilon\sum_{i=1}^kB_i\right)\otimes I,
\end{split}
\end{equation*}
hence $\Re L_{B}(X)$ is lower bounded on
\begin{equation*}
\mathcal{N}^{\perp}\left(\sum_{i=0}^k B_i\right)\otimes E=\left(\mathcal{N}^{\perp}(B_0)\cup \mathcal{N}^{\perp}\left(\sum_{i=1}^k B_i\right)\right)\otimes E,
\end{equation*}
thus $\Re L_{B}(X)$ is invertible on $\mathcal{N}^{\perp}(\sum_{i=0}^k B_i)\otimes E$. Since $\Im L_{B}(X)=\sum_{i=1}^kB_i\otimes \Im X_i$ it follows that $\mathcal{N}^{\perp}(\Im L_{B}(X))\leq \mathcal{N}^{\perp}(\Re L_{B}(X))$. Thus
\begin{equation*}
\begin{split}
L_{B}(X)&=\Re L_{B}(X)+i\Im L_{B}(X)\\
&=(\Re L_{B}(X))^{1/2}[I+i(\Re L_{B}(X))^{-1/2}\Im L_{B}(X)(\Re L_{B}(X))^{-1/2}](\Re L_{B}(X))^{1/2}
\end{split}
\end{equation*}
is lower bounded on $\mathcal{N}^{\perp}(\sum_{i=0}^k B_i)\otimes E$, since for any $A\in\mathcal{B}(\mathcal{N}^{\perp}(\Im L_{B}(X))\otimes E)$ we have $(I+iA)^*(I+iA)=I+A^*A$ which is lower bounded. It follows that $P^{\perp}L_{B}(X)P^{\perp}$ is lower bounded as well on $P^{\perp}\left(\mathcal{N}^{\perp}(\sum_{i=0}^k B_i)\otimes E\right)$, thus invertible and this implies that $\mathcal{S}(L_{B}(X))$ exists and is bounded. Since $X\in\Sigma^k$ and $\Re X_i\geq \epsilon I$, there exists an $\alpha\in[0,\pi/2)$ such that $W(X_i)\subseteq S_{\alpha}$ for all $1\leq i\leq k$. Notice that $W(I)\subseteq S_{\alpha}$ as well. Then by Lemma~\ref{L:sectorial} we have that $W(L_{B}(X))\subseteq S_{\alpha}$ and Proposition~\ref{P:Drury} implies \eqref{eq:P:DruryPsi:1}.

Now assume that $X\in\Pi^k$. Then there exists an $\epsilon>0$ such that $\Im X_i\geq \epsilon I$. Then there exists an $c_1>0$ such that
\begin{equation*}
\begin{split}
\Im L_{B}(X)&=\sum_{i=1}^kB_i\otimes \Im X_i\\
&\geq \sum_{i=1}^kB_i\otimes \epsilon I_{\mathcal{N}^{\perp}(\sum_{i=1}^k B_i)}\geq c_1 I_{\mathcal{N}^{\perp}(\sum_{i=1}^k B_i)}\otimes I,
\end{split}
\end{equation*}
hence $\Im L_{B}(X)$ is lower bounded on $\mathcal{N}^{\perp}\left(\sum_{i=1}^k B_i\right)\otimes E$. There also exists a $c_2\in\mathbb{R}$ such that $\Re L_{B}(X)-B_0\otimes I=\sum_{i=1}^kB_i\otimes \Re X_i\geq c_2 I_{\mathcal{N}^{\perp}(\sum_{i=1}^k B_i)}\otimes I$. Then there exists a $\theta\in(-\pi/2,0]$ such that
\begin{equation}\label{eq:P:DruryPsi:2}
\Re e^{i\theta}\Re L_{B}(X)-\Im e^{i\theta}\Im L_{B}(X)>0,
\end{equation}
equivalently lower bounded, on $(\mathcal{N}^{\perp}(B_0)\cup \mathcal{N}^{\perp}(\sum_{i=1}^k B_i))\otimes E$. Indeed, since $\Re L_{B}(X)\geq B_0\otimes I+c_2 I_{\mathcal{N}^{\perp}(\sum_{i=1}^k B_i)}\otimes I$ and $\Im L_{B}(X)\geq c_1 I_{\mathcal{N}^{\perp}(\sum_{i=1}^k B_i)}\otimes I$ and $B_0\geq 0$ is lower bounded on $\mathcal{N}^{\perp}(B_0)$, it suffices to have
\begin{equation*}
c_2\Re e^{i\theta}-c_1\Im e^{i\theta}>0.
\end{equation*}
Then $e^{i\theta} I$ and each $e^{i\theta} X_i$ for all $1\leq i\leq k$ is sectorial for some $\alpha\in[0,\pi/2)$, i.e. $W(e^{i\theta} I),W(e^{i\theta} X_i)\subseteq S_{\alpha}$. Thus Lemma~\ref{L:sectorial} implies that $W(e^{i\theta} L_{B}(X))\subseteq S_{\alpha}$ and \eqref{eq:P:DruryPsi:2} implies that $\Re (e^{i\theta} L_{B}(X))$ is lower bounded. Now a similar argument as in the previous case implies that $P^{\perp}e^{i\theta} L_{B}(X)P^{\perp}$ is lower bounded on $P^{\perp}\left((\mathcal{N}^{\perp}(B_0)\cup \mathcal{N}^{\perp}(\sum_{i=1}^k B_i))\otimes E\right)$ as well, thus invertible and this implies that $\mathcal{S}(e^{i\theta} L_{B}(X))$ exists and is bounded, which implies the existence and boundedness of $\mathcal{S}(L_{B}(X))$ as well, since $\mathcal{S}(e^{i\theta} L_{B}(X))=e^{i\theta} \mathcal{S}(L_{B}(X))$. Then Proposition~\ref{P:Drury} implies \eqref{eq:P:DruryPsi:1} with $e^{i\theta} L_{B}(X)$ in place of $L_{B}(X)$. Then \eqref{eq:P:DruryPsi:1} follows as well.

The above argumentation excludes the case when $S$ and $\mathcal{N}^{\perp}(B_0)\cup \mathcal{N}^{\perp}(\sum_{i=1}^k B_i)$ have trivial intersection. However in this case there is nothing to prove since then we have $\mathcal{S}(L_{B}(X))=L_{B}(X)$ and \eqref{eq:P:DruryPsi:1} is satisfied with $\alpha=0$.
\end{proof}

Now we may perform the free analytic continuation of our free function $F(X)=(w\otimes I)(\Psi_{F}(X))$ given in Theorem~\ref{T:Loewner_formula}.

\begin{proposition}\label{P:analyitic_cont}
Let $F:\mathbb{P}^k\mapsto \mathbb{P}$ be an operator monotone function and let $c_2>c_1>0$ be real numbers. Let $w\in\mathcal{B}_1^+(\mathcal{H})$ be a state in Theorem~\ref{T:Loewner_formula} such that for all $\dim(E)<\infty$ and $X\in\mathbb{P}_{c_1,c_2}(E)^k$ we have
\begin{equation}\label{eq:P:analyitic_cont:0}
F(X)=(w\otimes I)(\Psi_{F}(X)).
\end{equation}
Then $F$ has a free analytic continuation $\tilde{F}$ to the whole of $\Sigma^k$ and $\Pi^k$ for any Hilbert space $E$, such that $F$ maps $\Pi^k$ to $\Pi$, moreover $\tilde{F}$ is operator monotone and operator concave on $\mathbb{P}(E)^k$ and maps $\mathbb{P}(E)^k$ to $\mathbb{P}(E)$.
\end{proposition}
\begin{proof}
First we take care of the analytic continuation to $\Sigma^k$. Notice that according to Theorem~\ref{T:reconstruct_schur} and Proposition~\ref{P:separating_pencil} each direct summand of $\Psi_{F}(X)$ is the Schur complement of a linear pencil of the form
\begin{equation*}
L_{F,A,v}(X):=B(F,A,v)_0\otimes I+\sum_{i=1}^kB(F,A,v)_i\otimes (X_i-I).
\end{equation*}
For each $L_{F,A,v}$ we have that $B(F,A,v)_i\in\mathcal{B}^+(N_{F,A,v})_{*}$ for some finite dimensional Hilbert subspace $N_{F,A,v}\leq\mathcal{H}$, $\sum_{i=1}^kB(F,A,v)_i\leq B(F,A,v)_0$ and $\tr\{B(F,A,v)_0\}\leq \frac{F(c_2,\ldots,c_2)}{\min(1,c_1)}$ by Proposition~\ref{P:separating_pencil}. For a fixed $X\in\Sigma^k$ or $X\in\Pi^k$ it is easy to find a uniform upper bound on the norms of all $L_{F,A,v}(X)$ since we have the bound $\tr\{B(F,A,v)_0\}\leq \frac{F(c_2,\ldots,c_2)}{\min(1,c_1)}$ and $B(F,A,v)_i\in\mathcal{B}^+(\mathcal{H})_{*}$, $\sum_{i=1}^kB(F,A,v)_i\leq B(F,A,v)_0$. Then for all $A$ and $v$ separately, Proposition~\ref{P:DruryPsi} applies for $L_{F,A,v}(X)$ with subspace $S$ spanned by $v$. Then according to Definition~\ref{D:Psi} $\Psi_{F}(X)$ is the direct sum of all the Schur complements $\mathcal{S}(L_{F,A,v}(X))$, each with norm bounded from above by a uniform constant. Hence
\begin{equation*}
\Psi_{F}(X)^*\Psi_{F}(X)\leq K\sec^4\alpha I\otimes I
\end{equation*}
for some large enough real constant $K>0$ depending on $X$. Notice that $(w\otimes I):\mathcal{B}(\mathcal{H})\otimes \mathcal{B}(E)\mapsto \mathbb{C}\otimes \mathcal{B}(E)\cong \mathcal{B}(E)$ is a completely positive unital linear map. Hence by the Schwarz inequality for $2$-positive unital linear maps, see Proposition 3.3 \cite{paulsen}, we have
\begin{equation*}
\begin{split}
(w\otimes I)(\Psi_{F}(X))^*(w\otimes I)(\Psi_{F}(X))&\leq (w\otimes I)(\Psi_{F}(X)^*\Psi_{F}(X))\\
&\leq K\sec^4\alpha I,
\end{split}
\end{equation*}
hence \eqref{eq:P:analyitic_cont:0} defines a free holomorphic/analytic function on $\Sigma^k$ and $\Pi^k$.
 Then by Proposition~\ref{P:Schur_imgainary_invariant} $\tilde{F}$ maps $\Pi^k$ to $\Pi$ and by Proposition~\ref{P:Schur_concave_monotone} $\tilde{F}$ is operator monotone and operator concave on $\mathbb{P}(E)^k$ and clearly maps $\mathbb{P}(E)^k$ to $\mathbb{P}(E)$.
\end{proof}

\begin{theorem}\label{T:Loewner_several_var}
Let $E$ be a Hilbert space and let $F:\mathbb{P}(E)^k\mapsto \mathbb{P}(E)$ be a free function. Then the following are equivalent:
\begin{itemize}
\item[(a)]$F$ is operator monotone;
\item[(b)]$F$ is operator concave;
\item[(c)]There exists a Hilbert space $\mathcal{K}$, a closed subspace $\mathcal{K}_0\leq \mathcal{K}$ and the corresponding orthogonal projection $P_{\mathcal{K}_0}$ with range $\mathcal{K}_0$, $B_i\in\hat{\mathbb{P}}(\mathcal{K})$, $0\leq i\leq k$ with $B_0\geq \sum_{i=1}^kB_i$ and a state $w\in\mathcal{B}_1^+(\mathcal{K})^*$ such that for all $X\in\mathbb{P}(E)^k$ we have
\begin{equation}\label{eq:T:Loewner_several_var:0.1}
\begin{split}
F(X)=&w(B_{0,11})\otimes I+\sum_{i=1}^kw(B_{i,11})\otimes (X_i-I)\\
&-(w\otimes I)\left\{\left[B_{0,12}\otimes I+\sum_{i=1}^kB_{i,12}\otimes (X_i-I)\right]\right.\\
&\left[B_{0,22}\otimes I+\sum_{i=1}^kB_{i,22}\otimes (X_i-I)\right]^{-1}\\
&\left.\left[B_{0,21}\otimes I+\sum_{i=1}^kB_{i,21}\otimes (X_i-I)\right]\right\}
\end{split}
\end{equation}
where
\begin{equation}\label{eq:T:Loewner_several_var:0.2}
\begin{split}
{B}_{i,11}(A,v):=&P_{\mathcal{K}_0}{B}_{i}P_{\mathcal{K}_0},\\
{B}_{i,12}(A,v):=&P_{\mathcal{K}_0}{B}_{i}(I-P_{\mathcal{K}_0}),\\
{B}_{i,21}(A,v):=&(I-P_{\mathcal{K}_0}){B}_{i}P_{\mathcal{K}_0},\\
{B}_{i,22}(A,v):=&(I-P_{\mathcal{K}_0}){B}_{i}(I-P_{\mathcal{K}_0});
\end{split}
\end{equation}
\item[(d)]$F$ has a free analytic continuation to $\Pi^k$, mapping $\Pi^k$ to $\Pi$.
\end{itemize}
\end{theorem}
\begin{proof}
The equivalence between (a) and (b) is given by Theorem~\ref{T:concavemonotone} and Corollary~\ref{P:concavemonotone}.

That (c) implies the others is given by Proposition~\ref{P:Schur_imgainary_invariant}, Proposition~\ref{P:Schur_concave_monotone} and by Proposition~\ref{P:analyitic_cont}.

The proof that (d) implies (a) can be found in \cite{pascoe} as Lemma 4.8, we provide here the proof for completeness. First of all, similarly to the one variable case, a differentiable free function $F$ is operator monotone if and only if its Fr\'{e}chet-derivative $DF(X)(H)\geq 0$ for tuples $X$ in the domain of $F$ and any tuple $H\geq 0$. Indeed if $F$ is operator monotone then
\begin{equation*}
DF(X)(H)=\lim_{t\to 0+}\frac{F(X+tH)-F(X)}{t}
\end{equation*}
and by monotonicity $F(X+tH)\geq F(X)$ for $H\geq 0$ and $t\geq 0$, hence $DF(X)(H)\geq 0$. Conversely
\begin{equation*}
F(X)-F(Y)=\int_0^1 DF((1-t)Y+tX)(X-Y)dt,
\end{equation*}
hence if $X\geq Y$ and $DF(Z)(H)\geq 0$ for any $H\geq 0$ and $Z$ in the domain of $F$, then the integral on the right is positive semidefinite, hence $F(X)\geq F(Y)$ as well, establishing that $F$ is operator monotone. To show that (d) implies (a) assume that (d) holds, but $F$ is not operator monotone. Then by the above $DF(X)(H)$ is not positive semidefinite for some $X\in\mathbb{P}(E)^k$ and $H\geq 0$, $H\in\mathbb{S}(E)^k$. Since
\begin{equation*}
\begin{split}
F(X+itH)&=F(X)+itDF(X)(H)+O(t^2),\\
\Im F(X+itH)&=tDF(X)(H)+O(t^2),
\end{split}
\end{equation*}
for this $X$ and $H$ it follows that $\Im F(X+itH)$ is not positive semidefinite for small enough $t\geq 0$, a contradiction.

Lastly we show that (a) implies (c). For any fixed $c_2>c_1>0$, by Theorem~\ref{T:Loewner_formula} (c) holds for all $\dim(E)<\infty$ and $X\in\mathbb{P}_{c_1,c_2}(E)^k$. We use Proposition~\ref{P:analyitic_cont} to free analytically continue $F$ to $\Sigma^k$. We can follow this procedure for any $c_2>c_1>0$. Now we have identity and uniqueness theorems for non-commutative power series expansions for free analytic functions, see Theorem 7.2, 7.8 and 7.9 in \cite{verbovetskyi}. The non-commutative power series expansion around an arbitrary point $A\in\Sigma^k$ is called a TT series in \cite{verbovetskyi} and the coefficients of TT series are uniquely determined by the directional derivatives of the free function at $A$, and the series uniformly converges on some open ball around $A$. Therefore for each $c_2>c_1>0$ the free analytic function of the form as in \eqref{eq:T:Loewner_several_var:0.1} has uniquely determined free analytic continuation to the connected set $\Sigma^k$ by expanding the Schur complement into a non-commutative power series converging uniformly on open balls around arbitrary points in $A\in\Sigma^k$ representing the functions in \eqref{eq:T:Loewner_several_var:0.1}. Hence for any $c_2>c_1>0$ the formula given in (c) that we obtain by Theorem~\ref{T:Loewner_formula} all agree on $\Sigma^k$, and by similar arguments on $\Pi^k$ as well. Thus we can pick any one them of the form as in (c) for some fixed $c_2>c_1>0$ and for all $\dim(E)<\infty$. Now assume that $E$ is separable. Then by Lemma~\ref{L:strong_continuity} $F$ is strong operator continuous, hence for any bounded from above increasing net of $k$-tuple of operators $\{A_i\}_{i\in \mathcal{I}}$ with $A_i\in\mathbb{P}(E)^k$ we have 
\begin{equation}\label{eq:T:Loewner_several_var:1}
\sup_{i\in\mathcal{I}}F(A_i)=F\left(\sup_{i\in\mathcal{I}}A_i\right).
\end{equation}
If $E$ is non-separable, then by Assumption~\ref{A:semi-continuity2} for any bounded from above increasing net of $k$-tuple of operators $\{A_i\}_{i\in \mathcal{I}}$ with $A_i\in\mathbb{P}(E)^k$ we have 
\begin{equation*}
\sup_{i\in\mathcal{I}}F(A_i)\geq F\left(\sup_{i\in\mathcal{I}}A_i\right),
\end{equation*}
which implies \eqref{eq:T:Loewner_several_var:1}, since $F(A_i)\leq F\left(\sup_{i\in\mathcal{I}}A_i\right)$ by operator monotonicity of $F$ alone. Now using this strong continuity we can approximate any tuple $X\in\mathbb{P}(E)^k$ in the strong product topology by a bounded from above increasing net of tuples of finite rank operators $X_\gamma\in\mathbb{P}(E)^k$, such that $X_\gamma\to X$ strongly. For each $F(X_\gamma)$ (c) holds if we restrict each $X_\gamma$ to the orthogonal complement of its kernel. Notice that a free function defined by the formula \eqref{eq:T:Loewner_several_var:0.1} is strong operator continuous on norm bounded sets since it can be expanded into a uniformly norm-convergent non-commutative power series. The other way to see this is by Lemma~\ref{L:strong_continuity} and the norm continuity of the formula \eqref{eq:T:Loewner_several_var:0.1}. Then by the strong operator continuity of $F$, we have $F(X_\gamma)\to F(X)$ in the strong operator topology, hence (c) holds for any $E$ and $X\in\mathbb{P}(E)^k$.
\end{proof}

Now given an arbitrary operator monotone NC function $F:\mathbb{S}_{a_1,b_1}(E)\times\cdots\times\mathbb{S}_{a_k,b_k}\mapsto\mathbb{S}$ where $\mathbb{S}_{a,b}:=\{X\in\mathbb{S}(E):aI\leq X_i\leq bI\}$ for real numbers $b>a$, we can apply Theorem~\ref{T:Loewner_several_var} to free analytically continue the function using the following transformations. M\"obius transformations
\begin{equation}
g(x)=\frac{ax+b}{cx+d}
\end{equation}
are operator monotone for $x\in\mathbb{R}, x\neq -d/c$ provided $ad-bc>0$ and they map the upper complex half-plane into itself, see for example \cite{bendatsherman}. It is a relatively simple claculation to show that they also map $\Pi$ to $\Pi$. We can find a M\"obius transformation $g$ which maps $(a,b)$ to $(0,\infty)$ bijectively, with $g$ (and also its inverse) being operator monotone. Hence we may compose with such operator monotone transformations for each variable of $F$, to obtain a new operator monotone NC function $\hat{F}:\mathbb{P}^k\mapsto\mathbb{S}$ from $F$. If $\hat{F}$ is bounded from below then after adding some $cI$ for real $c\geq 0$ we get that $\hat{F}+cI>0$ for which we can apply Theorem~\ref{T:Loewner_several_var} to obtain the free analytic continuation mapping $\Pi^k$ to $\Pi$. If $\hat{F}$ is not bounded from below on $\mathbb{P}^k$, then $\overline{F}(X):=\hat{F}(X+\epsilon I)$ is bounded from below on $\mathbb{P}^k$ by operator monotonicity for any real $\epsilon>0$, hence $\overline{F}$ satisfies the previous case. In this way we obtain the free analytic continuation of the original $F$ mapping $\Pi^k$ to $\Pi$.

\section*{Acknowledgment}
An initial version of this paper appeared online in May 2014, where the author claimed more or less similar results stated in this paper. However the proofs contained a major flaw. The flaw was found later by the author, however theoretical counterexamples were communicated to the author earlier by James E. Pascoe and certain numerical counterexamples by Takeaki Yamazaki. The author would like to express his gratitude to James E. Pascoe and Takeaki Yamazaki for providing counterexamples. The techniques appearing in this current paper are very different from those in the initial one.

This work was partly supported by the Research Fellowship of the Canon Foundation, SGU project of Kyoto University, the JSPS international research fellowship grant No. 14F04320, the "Lend\"ulet" Program (LP2012-46/2012) of the Hungarian Academy of Sciences and the National Research Foundation of Korea (NRF) grant funded by the Korea government (MEST) No. 2015R1A3A2031159.




\end{document}